\documentclass[12pt,a4paper]{article}
\title{A Hybrid-High Order Method for Quasilinear Elliptic Problems of Nonmonotone Type}

\author{Thirupathi Gudi\footnote{Department of Mathematics, Indian Institute of Science, Bangalore 560012, India. Email. gudi@iisc.ac.in},
\; Gouranga Mallik\footnote{Department of Mathematics, Indian Institute of Science, Bangalore 560012, India. Email. gourangam@iisc.ac.in}
\; and Tamal Pramanick\footnote{Department of Mathematics, National Institute of Technology Calicut, Kozhikode 673601, India. Email. tamal@nitc.ac.in}}

 \usepackage{amsmath,amsthm,amssymb,enumerate}
\usepackage{gensymb}
\usepackage[square,sort&compress,comma,numbers]{natbib}

\usepackage[sc]{mathpazo}
\usepackage{multirow} 

\usepackage{multicol}

\usepackage{newtxtext,newtxmath}
\usepackage{subfig}
\usepackage{graphicx}
\usepackage{epstopdf}
\usepackage{hyperref}
\usepackage{cancel}
\usepackage[margin=3cm]{geometry}
\usepackage{tikz}
\usepackage{caption}
\usetikzlibrary{shapes,calc}
\usepackage{verbatim}
\usepackage{mathrsfs}
\usepackage{algorithm}
\usepackage{accents}

\chardef\bslash=`\\ 

\newtheorem{theorem}{Theorem}[section]
\newtheorem{lemma}[theorem]{Lemma}
\newtheorem{proof of lemma}[theorem]{Proof of Lemma}

\theoremstyle{definition}

\newtheorem{remark}[theorem]{Remark}
\numberwithin{equation}{section}

\newcommand{\dx}{{\rm\,dx}}

\newcommand{\ds}{{\rm\,ds}}
\newcommand{\dt}{{\rm\,dt}}

\newcommand{\fl}{\quad\forall}

\newcommand{\bR}{\mathbb R}

\newcommand{\mbP}{\mathbb P}

\newcommand{\cF}{\mathcal F}

\newcommand{\cT}{\mathcal T}

\newcommand{\cN}{\mathcal N}

\newcommand{\Gtk}{\boldsymbol{G}_T^{k}}

\newcommand{\Utk}{\underline{U}_T^k}
\newcommand{\Ghk}{\boldsymbol{G}_h^{k}}

\newcommand{\lt}{L^2(\Omega)}

\newcommand{\Holders}{H\"{o}lder's }

\newcommand{\GardingsType}{G$\mathring{\text{a}}$rding-type }
\newcommand{\CS}{Cauchy--Schwarz }

\newcommand{\integ}{\int_\Omega}
\newcommand{\bfnTF}{\boldsymbol{n}_{TF}}
\newcommand{\bfn}{\boldsymbol{n}}
\newcommand{\sit}{\sum_{T\in\mathcal{T}_h}\int_T}

\newcommand{\sumt}{\sum_{T\in\mathcal{T}_h}}

\newcommand{\Ctr}{C_{\text{tr}}}
\newcommand{\Ctrc}{C_{\text{tr},\text{c}}}

\allowdisplaybreaks

\begin{document}
\date{\today}
\maketitle

\begin{abstract} 
In this paper, we design and analyze a Hybrid-High Order (HHO) approximation for a class of quasilinear elliptic problems of nonmonotone type. The proposed method has several advantages, for instance, it supports arbitrary order of approximation and general polytopal meshes. The key ingredients involve local reconstruction and high-order stabilization terms. Existence and uniqueness of the discrete solution are shown by Brouwer's fixed point theorem and contraction result. A priori error estimate is shown in discrete energy norm that shows optimal order convergence rate. Numerical experiments are performed to substantiate the theoretical results.
\end{abstract}

\section{Introduction}\label{sec:intro}
We consider here numerical approximation for the nonlinear elliptic boundary value problem:
\begin{subequations}\label{quasilin_pde_intro}
	\begin{align}
		-\nabla{\cdot}(a(x,u)\nabla u)&= f(x)\quad\text{in}~~ \Omega,\\
		u(x)&=0\quad\text{on}~~ \partial\Omega,
	\end{align}
\end{subequations}
where $\Omega$ is a bounded convex polytopal domain in $\mathbb{R}^d$, $d\in \{2,3\}$ with Lipschitz boundary $\partial\Omega$ and $a:\bar{\Omega}\times\bR \to \bR$ is a nonlinear function of its arguments.  For simplicity, the homogeneous boundary condition is considered. Some additional  assumptions are stated in the appropriate section.  The main purpose of the article is to devise and analyse Hybrid-High Order (HHO) approximation for the problem~\eqref{quasilin_pde_intro} on general meshes inspired by the HHO methods of \cite{Piet_Ern_Lem_14_arb_local} for linear diffusion model problem, \cite{Piet_Dron_Ern_15_HHO_Adv_Diff_Rea} for degenerate advection–diffusion–reaction models and \cite{Piet_Dron_16_Leray} for nonlinear steady Leray–Lions equation.
The proposed method became very famous over the last decade. The method consists several advantages, such as the HHO discretization supports general polytopal meshes, and it allows arbitrary order of approximation of polynomials. The HHO method complies with physics and it is robust with respect to the variations of physical coefficients. It focuses on reproduction at the discrete level of key continuous properties such as local balances and flux continuity. The computational cost for HHO method can be reduced by the use of compact stencil and static condensation.

Over the last several years, a major focus has been on analyzing discretization methods for partial differential equations (PDEs) which support arbitrary-order discretization on general meshes including nonmatching interfaces and polytopal cells. A lot of research has been done on arbitrary-order polytopal method for linear diffusion equation. To name a few, we refer \cite{Anton_Gian_Hous_13_DG_Complicated_Dom,Bassi_Bott_Colo_Piet_12_DG_Physical} for the constructions of polyhedral elements for adaptive mesh coarsening, Hybridizable Discontinuous Galerkin method of \cite{Cock_Gopal_Laz_09_Unif,Cast_Cock_Peru_Scho_00_local_DG}, the Virtual Element method of \cite{Beira_Brezz_Cang_Man_Mar_13_VEM,Beira_Brezz_Mari_13_VEM_Elast,Brezz_Falk_Mari_14_VEM_Mixed}, the High-Order Mimetic method of \cite{Lip_Manz_14_Mimetic}, the Weak Galerkin method of \cite{Mu_Wang_Ye_15_WeakGaler_Nonlin,Wang_Ye_13_Weak_Galerkin,Wang_Ye_14_Mixed_Weak_Galerkin}, the Gradient Discretisation methods of \cite{Dron_Eym_Herb_GDM_16,Dron_Eym_GDM_18_book,Piet_Dron_Man_18_GradDis_Polytope} and the Multiscale Hybrid-Mixed method of \cite{Arya_Hard_Pare_Vale_13_Multiscale}. Recently, robust HHO approximation scheme for Poisson problem on polytopal meshes with small edges/faces is analyzed in \cite{Jerome_Liam_21_hTscaling}, where the standard scaling of $h_F$ for the stabilization is replaced by $h_T$.
There are some connections of HHO method with the Hybridizable Discontinuous Galerkin (HDG) method but the choice of stabilization for HHO method is different from that of HDG as to deliver higher-order convergence rate for HHO discretization.
On the other side, if we consider the nonconforming Virtual Element Methods (ncVEM), the devising viewpoint with HHO is different. Nonconforming VEM considers the computable projection of virtual functions instead of a reconstruction operator, and the stabilization for both methods achieves similar convergence rates as HHO but is written differently. The close connections of HHO method with HDG \cite{Cock_Gopal_Laz_09_Unif} and ncVEM \cite{Agus_Lip_Manz_16_NCVEM}  have been analyzed in \cite{Cock_Piet_Ern_16_HHO_HDG}.

HHO method in lowest-order case belongs to the Hybrid Mixed Mimetic family \cite{Dron_Eym_Gall_Her_10_Mimetic}, which includes the mixed-hybrid Mimetic Finite Differences \cite{Brez_Lipn_Simo_05_Mim_Dif}, the Hybrid Finite Volume \cite{Eym_Gall_Herb_10_NonconfMesh} and the Mixed Finite Volume \cite{Dron_14_FV_Diff, Dron_Eym_06_MixedFV}. Recently, the HHO method has been bridged in \cite{Lem_21_Bridging_HHO_VEM} with virtual element method. Some more related approaches can also be found in \cite{Bone_Ern_14_poly, Dron_Eym_Gall_10_unified,Brezzi_Lip_Shas_05_Conv_Mimetic, Brezzi_Lip_Shas_07_polyhedral, Kuzn_Lip_Shas_04_Mimetic_poly}.
Several works for HHO methods involving linear and nonlinear PDEs can be found in various articles, such as pure diffusion \cite{Piet_Ern_Lem_14_arb_local}, advection-diffusion \cite{Piet_Dron_Ern_15_HHO_Adv_Diff_Rea}, interface problems \cite{Burman_Ern_18_HHO_Interface}, viscosity-dependent Stokes problem \cite{Piet_Ern_Link_16_HHO_Stokes} for linear PDEs, elliptic obstacle problem \cite{Cicu_Ern_Gudi_20_HHO_Obstacle},
a nonlinear elasticity with infinitesimal deformations \cite{Botti_Piet_Socha_17_HHO_Nonlin_Elast}, steady incompressible Navier Stokes equations \cite{Piet_Krell_18_NSE} and Leray-Lions operators \cite{Piet_Dron_16_Leray,Piet_Dron_Harn_21_Improved_Leray} for nonlinear PDEs.

The quasilinear problem of nonmonotone type \eqref{quasilin_pde_intro} can be thought of as the stationary heat problem with variable nonlinear diffusion coefficient. This has many engineering applications, for instance, heat distribution for metal bodies.  Finite element approximations for the nonlinear problem \eqref{quasilin_pde_intro} has been studied in \cite{Gudi_AKP_07_DG_quasi,Gudi_NN_AKP_08_HpDG_Quasi} for a priori error estimate using DG and $hp$-DG methods, and see also  \cite{Chun_Vict_13_Apost_DG_Quasi,Chun_Vict_09_APost_FV_Quasi,Chun_Ming_12_DGFV_Quasi,Song_Zhang_15_SupConv} for various a priori and a posteriori error estimate.  For some related works on $hp$-DG method for strongly nonlinear elliptic problem, we refer \cite{Gudi_NN_AKP_08_Strongly_nonlin,Chun_Wang_Lin_18_Two_Grid} and references therein. We also refer articles on various nonlinear problem on second-order elliptic PDEs, Weak Galerkin methods of \cite{Mu_Wang_Ye_15_WeakGaler_Nonlin}, Mimetic finite difference approximation of \cite{Anto_Big_Ver_15_Mimetic_quasilin} and Virtual element method of \cite{Canggiani_20_VEM_Quasilin}. The finite element approximations for the quasilinear problem are studied under various regularity assumptions on the coefficient $a$ and on the solution $u$, for instance, $a$ in $C^2$ and $u\in W^{2,2+\epsilon}(\Omega), \epsilon>0$ or $u\in H^2(\Omega)$, see \cite{Xu_96_TwoGrid_Nonlin,Gudi_AKP_07_DG_quasi,Chun_Vict_11_TwoGridDG_Quasi,Chun_Vict_13_Apost_DG_Quasi,Dou_Dup_Ser_71_DivExitUniq,Lip_Michal_96,Milner_85}.

In this article, we establish optimal order a priori error estimate in discrete energy norm for the hybrid high-order approximation for quasilinear elliptic problem of nonmonotone type. We assume the solution $u\in H^1_0(\Omega)$ of \eqref{quasilin_pde_intro} belongs to $H^2(\Omega)\cap W^{1,\infty}(\Omega)$ for $d=2$ and belongs to $H^3(\Omega)$ for $d=3$. We use local reconstruction and high-order stabilization in the discrete formulation.  First, we establish existence, uniqueness and error estimate for an auxiliary second-order nonselfadjoint linear elliptic problem using \GardingsType inequality. The well-posedness of this auxiliary problem helps us to formulate a suitable nonlinear map which posesses a ball to ball mapping and contraction properties. For sufficiently small mesh parameter, the existence, uniqueness and error estimate for the solution of HHO approximation of quasilinear problem are shown by  Brouwer's fixed point theorem and contraction result on quasi-uniform mesh.

The organization of the paper is as follows. Section~\ref{sec:intro} is introductory in nature, and Section~\ref{sec:HHO_Dis} is devoted to notation, definition and preliminaries related to HHO discretization. In Section~\ref{sec:Nonself_Adj}, we have discussed the HHO approximation for a linear nonselfadjoint elliptic problem and established error estimates.
Section~\ref{sec:Quasilin} is devoted to HHO approximation for the solution of quasilinear elliptic problem \eqref{quasilin_pde}. In Section~\ref{Sec:Num_Quasilin}, numerical experiments are performed to illustrate the theoretical results obtained in this article. Finally, in Section~\ref{sec:Conclusion} we present summary of the article and describe some possible extensions. 

Throughout the paper, standard notation on Lebesgue and Sobolev spaces and their norms are employed.
The standard semi-norm and norm on $H^{s}(\Omega)$ (resp. $W^{s,p} (\Omega)$) for $s>0$ are denoted by $|\bullet|_{s}$ and $\|\bullet\|_{s}$ (resp. $|\bullet|_{s,p}$ and $\|\bullet\|_{s,p}$ ). The positive constants $C$ appearing in the inequalities denote generic constants which do not depend on the mesh-size. The notation $a\lesssim b$ means that there exists a generic constant $C$ independent of the mesh parameters such that $a \leq Cb$; $a\approx b$ abbreviates $a\lesssim b\lesssim a$.

\section{Hybrid-High Order discretization}\label{sec:HHO_Dis}

\subsection{Discrete setting}
We consider a sequence of refined meshes $(\cT_h)_{h>0}$ where the parameter $h$ denotes the meshsize and goes to zero during the refinement process. For all $h>0$, we assume that the mesh $\cT_h$ covers $\Omega$ exactly and consists of a finite collection of non-empty disjoint open polyhedral cells $T$ such that $\overline{\Omega}=\cup_{T\in\cT_h} \overline{T}$, and $h=\max_{h\in\cT_h} h_T$, where $h_T$ is the diameter of $T$. A closed subset $F$ of $\Omega$ is defined to be a mesh face if it is a subset of an affine hyperplane $H_F$ with positive $(d-1)$-dimensional Hausdorff measure and if either of the following two statements holds true: (i) There exist $T_1(F)$ and $T_2(F)$ in $\cT_h$ such that $F\subset\partial T_1(F)\cap \partial T_2(F)\cap H_F$; in this case, the face $F$ is called an internal face; (ii) There exists $T(F)\in\cT_h$ such that $F \subset \partial T(F) \cap \partial\Omega \cap H_F$; in this case, the face $F$ is called a boundary face. The set of mesh faces is a partition of the mesh skeleton, i.e., $\cup_{T\in\cT_h}\partial T=\cup_{F\in\cF_h}\bar{F}$, where $\cF_h:=\cF_h^i\cup \cF_h^b$ is the collection of all faces that is the union of the set of all the internal faces $\cF_h^i$ and the set of all the boundary faces $\cF_h^b$. Let $h_F$ denote the diameter of $F\in \cF_h$. For each $T\in \cT_h$, the set $F_T:=\{F\in \cF_h\, |\, F\subset\partial T\}$ denotes the collection of all faces contained in $\partial T$, $\bfn_T$ the unit outward normal to $T$, and we set $\bfnTF:=\bfn_T|_F$ for all $F\in \cF_h$. Following \cite[Definition 1]{Piet_Ern_15_mesh}, we assume that the mesh sequence $(\cT_h)_{h>0}$ is admissible, in the sense that, for all $h>0$, $\cT_h$ admits a matching simplicial submesh $\cT_h$ (i.e., every cell and face of $\cT_h$ is a subset of a cell and a face of $\cT_h$, respectively) so that the mesh sequence $(\cT_h)_{h>0}$ is shape-regular in the usual sense and all the cells and faces of $\cT_h$ have uniformly comparable diameter to the cell and face of $\cT_h$ to which they belong. Owing to \cite[Lemma~1.42]{Piet_Ern_12_DG_book}, for $T\in\mathcal{T}_h$ and $F\in\mathcal{F}_T$, $h_F$ is comparable to $h_T$ in the sense that
\begin{align*}
	\varrho^{2}h_T\leq h_F\leq h_T,
\end{align*}
where $\varrho$ is the mesh regularity parameter. Moreover, there exists an integer $N_{\partial}$ depending on $\varrho$ and $d$ such that (see \cite[Lemma~1.41]{Piet_Ern_12_DG_book})
\begin{align*}
	\underset{T\in\mathcal{T}_h}{\text{max}}\, \text{card} (\mathcal{F}_T)\leq N_{\partial}.
\end{align*}

\noindent There also exist real numbers $\Ctr$ and $\Ctrc$ depending on $\varrho$ but independent of $h$ such that the following discrete and continuous trace inequalities hold for all $T\in\cT_h$ and $F\in\cF_T$ (see \cite[Lemma~1.46 and 1.49]{Piet_Ern_12_DG_book})
\begin{align}
	\|v\|_F&\leq \Ctr h_F^{-1/2} \|v\|_T\quad\forall v\in\mathbb{P}_d^l(T),\label{dis_trace}\\ 
	\|v\|_{\partial T}&\leq C_{\text{tr},\text{c}}(h_T^{-1}\|v\|_T^2+ h_T\|\nabla v\|_T^2)^{1/2}\quad\forall v\in H^1(T),
\end{align}
where $\mathbb{P}_d^l(T)$ is the space of polynomial of degree at most $l$ on $T\in\cT_h$. There exists a real number $C_{\text{app}}$ depending on $\varrho$ and $l$ but independent of $h$ such that, for all $T\in\cT_h$, denoting by $\pi_T^l$ the $L^2$-orthogonal projector on $\mathbb{P}_d^l(T)$, the following holds (see \cite[Lemma~1.58 \& 1.59]{Piet_Ern_12_DG_book}): For all $ s\in \{1,\ldots,l+1\}$ and all $v\in H^s(T)$,
\begin{align}
	|v-\pi_T^l v|_{H^m(T)}+ h_T^{1/2}|v-\pi_T^l v|_{H^m(\partial T)}\leq C_{\text{app}} h_T^{s-m}|v|_{H^s(T)},\quad \forall m\in \{0, \ldots,s-1\}.\label{proj_est}
\end{align}

\subsection{Discrete spaces}
Let $k\geq 0$ be a fixed polynomial degree. For $T\in\cT_h$, define the local space of degrees of freedom (DOFs) by
\begin{align}\label{local_dofs}
	\underline{U}_T^k&:= \mathbb{P}_d^k(T)\times \left\{\underset{F\in\mathcal{F}_T}{\times} \mathbb{P}_{d-1}^k(F)\right\},
\end{align}
where $\mathbb{P}_d^k(T)$ is the space of polynomials of degree at most $k$ on $T\in\cT_h$ and 
$\mathbb{P}_{d-1}^k(F)$ is the space of polynomial of degree at most $k$ on the face $F\in\cF_h$. 
The global space of DOFs is obtained by patching interface values in \eqref{local_dofs} as
\begin{align*}
	\underline{U}_h^k&:= \left\{\underset{T\in\mathcal{T}_h}{\times} \mathbb{P}_d^k(T)\right\}\times \left\{\underset{F\in\mathcal{F}_h}{\times} \mathbb{P}_{d-1}^k(F)\right\}.
\end{align*}
The zero boundary condition can be imposed in the above discrete space $\underline{U}_h^k$ as follows:
\begin{align*}
	\underline{U}_{h,0}^k&:= \left\{\underline{v}_h=\left((v_T)_{T\in\cT_h},(v_F)_{F\in\cF_h}\right)\in\underline{U}_h^k  \,|\, v_F\equiv 0\fl F\in\mathcal{F}_h^b\right\}.
\end{align*}
For $\underline{v}_h\in \underline{U}_h^k$, we understand $v_h\in L^2(\Omega)$ by $v_h|_T=v_T$. 
The local interpolation operator $I_T^k:H^1(T)\to \underline{U}_T^k$, is such that, for all $v\in H^1(T)$
\begin{align}
	I_T^kv:= (\pi_T^k v, (\pi_F^kv)_{F\in\mathcal{F}_T}),\label{defn_interpolant}
\end{align}
where $\pi_F^k$ is the $L^2$-orthogonal projector on $\mathbb{P}_{d-1}^k(F)$. The corresponding global interpolation operator $I_h^k:H^1(\Omega)\to \underline{U}_h^k$ is such that, for all $v\in H^1(\Omega)$,
\begin{align*}
	I_h^kv:= ((\pi_T^kv)_{T\in\mathcal{T}_h}, (\pi_F^kv)_{F\in\mathcal{F}_h}).
\end{align*}
When applied to $H_0^1(\Omega)$, $I_h^k$ maps onto $\underline{U}_{h,0}^k$. 

Below, we state Lebesgue embedding result. For a proof, we refer to \cite[Lemma 5.1]{Piet_Dron_16_Leray}.
\begin{lemma}[Direct and reverse Lebesgue embeddings]\label{lem_Sob_Inv_Ineq} Let $\cT_h$ be a regular mesh with $T\in\cT_h$. Let $k\in\mathbb{N}$ and $q,m\in [1,\infty]$. Then
	\begin{equation}
		\|w\|_{L^q(T)}\approx |T|^{\frac{1}{q}-\frac{1}{m}}\|w\|_{L^m(T)}\fl w\in \mathbb{P}^k(T).
	\end{equation}
\end{lemma}

Define the Sobolev exponent $p^*$ of $p$ by
\begin{align*}
	p^*:=
	\begin{cases}
		\frac{dp}{d-p}\quad \text{ if } p<d,\\
		+\infty\quad \text{ if } p\geq d.
	\end{cases}
\end{align*}
The next lemma on discrete Sobolev embeddings is stated from \cite[Proposition 5.4]{Piet_Dron_16_Leray} and this is used to obtain various boundedness result in the following sections.
\begin{lemma}[Discrete Sobolev embeddings]\label{lem_dis_emb} Let $(\cT_h)_{h>0}$ be an admissible mesh sequence of $\Omega\subset \bR^d$. Let $1\leq q\leq p^*$ if $1\leq p<d$ and $1\leq q<\infty$ if $p\geq d$. Then, there exists $C$ only depending on $\Omega, \varrho, k,q$ and $p$ such that
	\begin{equation*}
		\|v_h\|_{L^q(\Omega)}\leq C\|\underline{v}_h\|_{1,p,h}\fl\underline{v}_h\in \underline{U}_{h,0}^k, 
	\end{equation*}
	where $\displaystyle\|\underline{v}_h\|_{1,p,h}:=\left(\sum_{T\in\cT_h}\|\underline{v}_T\|_{1,p,T}\right)^{1/p}$ with $$\displaystyle \|\underline{v}_T\|_{1,p,T}:=\left(\|\nabla v_T\|^p_{L^p(T)^d}+\sum_{F\in\cF_T}h_F^{1-p}\|v_F-v_T\|_{L^p(T)}^p\right)^{1/p}.$$
	In particular, 
	\begin{equation}
		\|v_h\|_{L^6(\Omega)}\leq C\|\underline{v}_h\|_{1,2,h}\fl\underline{v}_h\in \underline{U}_{h,0}^k.\label{dis_emb}    
	\end{equation}
\end{lemma}
We use the abbreviation $\|\bullet\|_{1,h}$ for $\|\bullet\|_{1,2,h}$ in the subsequent analysis.

\subsection{Local reconstructions and stabilization operators}
In this subsection, some essential ingredients related to HHO formulation are defined. For $T\in\cT_h$, we define the local reconstruction operator $R_T^{k+1}:\underline{U}_T^k\to \mathbb{P}_d^{k+1}(T)$ such that, for $\underline{v}_T=(v_T,(v_F)_{F\in\cF_T})$,
\begin{subequations}\label{recons_oper}
	\begin{align}
		(\nabla R_T^{k+1}\underline{v}_T,\nabla w)_T&=(\nabla v_T,\nabla w)_T+ \sum_{F\in\mathcal{F}_T}(v_F-v_T,\nabla w{\cdot}\bfnTF)_F,\label{hho3}\\
		\left(R_T^{k+1}\underline{v}_T,1\right)_T&= \left(v_T,1\right)_T,\label{hho9}
	\end{align}
\end{subequations}
where \eqref{hho3} is enforced for all $w\in \mathbb{P}_d^{k+1}(T)$.  
A global reconstruction operator $R_h^{k+1}:\underline{U}_h^k\to\mathbb{P}_d^{k+1}(\cT_h)$ is defined by $R_h^{k+1} \underline{v}_h|_T=R_T^{k+1} \underline{v}_T$.

Define the local gradient reconstruction $\Gtk:\Utk\to \mathbb{P}^k(T)^d$ such that, for all $\underline{v}_T\in\underline{U}^k_T$,
\begin{align}\label{Grad_recons}
	(\Gtk\underline{v}_T,\boldsymbol{\tau})_T=(\nabla v_T,\boldsymbol{\tau})_T+\sum_{F\in\cF_T}(v_F-v_T,\boldsymbol{\tau}{\cdot}\bfnTF)_F\fl\boldsymbol{\tau}\in\mathbb{P}^k(T)^d.
\end{align}
Moreover, from \cite[Lemma~4.10]{Piet_Jero_HHO_Book_20}, it holds
\begin{align}\label{Grad_recons_proj}
	(\Gtk\underline{v}_T,\boldsymbol{\tau})_T=(\nabla v_T,\boldsymbol{\tau})_T+\sum_{F\in\cF_T}(v_F-v_T,(\pi_T^k\boldsymbol{\tau}){\cdot}\bfnTF)_F\fl\boldsymbol{\tau}\in L^1(T)^d.
\end{align}
The relation between $\Gtk$ and $R_T^{k+1}$ is established by taking $\boldsymbol{\tau}=\nabla w$ with $w\in\mbP_d^{k+1}(T)$ in \eqref{recons_oper} and comparing with \eqref{Grad_recons} as
\begin{align}
	(\Gtk\underline{v}_T-\nabla R_T^{k+1}\underline{v}_T,\nabla w)_T =0\fl w\in\mbP^{k+1}(T). 
\end{align}
In other words, $\nabla R_T^{k+1}\underline{v}_T$ is the $L^2$-orthogonal projection of $\Gtk\underline{v}_T$ on $\nabla \mbP^{k+1}(T)\subset P^k(T)^d$ and $\|\nabla R_T^{k+1}\underline{v}_T\|_T\leq \|\Gtk\underline{v}_T\|_T$.

For $F\in\cF_T$, define the local stabilization operator $S_F^k:\underline{U}_T^k\to \mathbb{P}_{d-1}^{k}(F)$ by
\begin{align}\label{stab_defn}
	S_F^k\underline{v}_T:=\pi_F^k\left(v_F-v_T-\left(R_T^{k+1}\underline{v}_T-\pi_T^k R_T^{k+1}\underline{v}_T\right)\right).   
\end{align}

The next lemma follows from the property $R_T^{k+1}I_T^k v=\pi_T^{1,k+1}v$ for $v\in W^{1,1}(T)$, where $\pi_T^{1,k+1}$ is the elliptic projector, see \cite[Definition~1.39]{Piet_Jero_HHO_Book_20} and its approximation property \cite[Theorem~1.48]{Piet_Jero_HHO_Book_20} as
\begin{lemma}[Approximation properties of $R_T^{k+1}I_T^k$]\label{lem_apprx_recons}
	There exists a real number $C>0$, depending on $\varrho$ but independent of $h_T$ such that, for all $v\in H^{s+1}(T)$ for some  $s\in\{0,1,\ldots,k\}$,
	\begin{align}
		&\|v-R_T^{k+1}I_T^k v\|_T+ h_T^{1/2}\|v-R_T^{k+1}I_T^k v\|_{\partial T}+ h_T\|\nabla(v-R_T^{k+1}I_T^k v)\|_T \notag\\
		&\qquad+h_T^{3/2}\|\nabla(v-R_T^{k+1}I_T^k v)\|_{\partial T}\leq Ch_T^{s+1} \|u\|_{H^{s+1}(T)}.\label{eqn_recons_est}
	\end{align}
\end{lemma}
The property $\Gtk I_T^k v=\pi_T^k(\nabla v)$ for $v\in W^{1,1}(T)$ and the approximation property for $L^2$ projector $\pi_T^k$ lead to
\begin{lemma}[Approximation properties of $\Gtk I_T^k$]\label{lem_apprx_Gtk}
	\cite[Lemma~3.24]{Piet_Jero_HHO_Book_20} There exists a real number $C>0$, depending on $\varrho$ but independent of $h_T$ such that, for all $v\in H^{s+1}(T)$ for some $s\in\{0,1,\ldots,k\}$,
	\begin{align}
		\|\nabla v-\Gtk I_T^k v\|_{L^2(T)}\leq Ch_T^{s} |u|_{H^{s+1}(T)}.\label{eqn_recons_Gtk}
	\end{align}
\end{lemma}

\section{Nonselfadjoint linear elliptic problems and error estimate}\label{sec:Nonself_Adj}
To analyze the existence and uniqueness of the discrete solution of \eqref{quasilin_pde_intro}, we require a linearized problem which is essentially a nonselfadjoint elliptic problem (the explicitly form is described in Section~\ref{sec:Quasilin}). This motivates us to consider a general nonselfadjoint elliptic partial differential equation:
\begin{subequations}\label{cts_lin_nonselfadj}
	\begin{align}
		-\nabla{{\cdot}}(a(x)\nabla u)+ \vec{b}(x){\cdot}\nabla u+a_0(x)u&= p(x)\quad\text{in}~~ \Omega,\\
		u&=0\quad\text{on}~~ \partial\Omega.
	\end{align}
\end{subequations}
In this section, we consider the HHO approximation and establish error estimate for the discrete solution of the above problem ~\ref{cts_lin_nonselfadj}. The existence and uniqueness results will be used in the error analysis of HHO approximation for quasilinear elliptic problem~\eqref{quasilin_pde_intro}. We adopt the following assumptions on the above problem~\eqref{cts_lin_nonselfadj}. 

\noindent\textbf{Assumption A:}
\begin{enumerate}
	\item There exists $\alpha>0$ such that $\alpha\leq a(x)$ and $a_0(x)\geq 0$ for all $x\in\overline{\Omega}$.
	\item $a\in W_{\infty}^1(\Omega)$, $\vec{b}\in\boldsymbol{L}^{\infty}(\Omega)$ and  $a_0\in L^{\infty}(\Omega)$ with\\
	$M :=\text{max}\{\|a\|_{L^{\infty}(\Omega)}, \|\vec{b}\|_{\boldsymbol{L}^{\infty}(\Omega)}, \|a_0\|_{L^{\infty}(\Omega)}\}$.
	\item $p\in L^2(\Omega)$.
\end{enumerate}
Then, it is well known (see \cite[Theorem~3.12]{Ern_Guer_04_FEMs} and \cite[Lemma~9.17]{Gil_Trud_book_83}) that there exists a unique solution $u\in H^2(\Omega)\cap H^1_0(\Omega)$ to the problem (\ref{cts_lin_nonselfadj}) and it satisfies the following a priori bound
\begin{align}
	\|u\|_{H^2(\Omega)}\leq C\|p\|_{L^2(\Omega)}.
\end{align}

\subsection{HHO discretization of nonselfadjoint linear elliptic problem}
\noindent Define a discrete counterpart of $(a_T\nabla u,\nabla v)_T$, where $a_T(x)=a(x)|_T$, as the local discrete bilinear form 
\begin{align}
	B_T(\underline{u}_T,\underline{v}_T)&:= \left(a_T\Gtk\underline{u}_T, \Gtk\underline{v}_T\right)_T+ s_T(\underline{u}_T,\underline{v}_T)\notag\\
	&\qquad+\left(\vec{b}{\cdot}\nabla R_T^{k+1}(\underline{u}_T), v_T\right)_T+ (a_0u_T,v_T)_T, \label{hho24}
\end{align}
where the stabilization term $s_T(\underline{u}_T,\underline{v}_T)$ is defined by
\begin{equation}\label{stabilization}
	s_T(\underline{u}_T,\underline{v}_T):= \sum_{F\in\mathcal{F}_T}\frac{a_{TF}^{\infty}}{h_F} \left(S_F^k\underline{u}_T,S_F^k\underline{v}_T\right)_F \text{ with } a_{TF}^{\infty}:=\|a_T\|_{L^{\infty}(F)}.
\end{equation}
The discrete formulation of \eqref{cts_lin_nonselfadj} seeks $\underline{u}_h\in \underline{U}_{h,0}^k$ such that
\begin{align}
	B_h(\underline{u}_h,\underline{v}_h)= (p,v_h)\quad\forall \underline{v}_h\in \underline{U}_{h,0}^k,\label{HHO_nonself}
\end{align}
where $B_h(\underline{u}_h,\underline{v}_h):= \sum_{T\in\mathcal{T}_h} B_T(\underline{u}_T,\underline{v}_T)$.

Define the energy seminorms on $\underline{U}_{h}^k$ (norms on  $\underline{U}_{h,0}^k$ owing to the homogeneous boundary condition) as follows:
\begin{align}
	&\|\underline{v}_h\|_{a,h}^2:= \sum_{T\in\mathcal{T}_h}\|\underline{v}_T\|_{a,T}^2\text{ and } \|\underline{v}_h\|_{1,h}:=\|\underline{v}_h\|_{1,2,h}\quad \text{ for all } \underline{v}_h\in \underline{U}_{h}^k,
\end{align}
where the local contributions are defined as
\begin{align}
	\|\underline{v}_T\|_{a,T}^2&:= \|a_T^{1/2}\Gtk \underline{v}_T\|_T^2+ \sum_{F\in\mathcal{F}_T}\frac{a_{TF}^{\infty}}{h_F} \|v_F-v_T\|_F^2.\label{norm_aT}
\end{align}
With the help of Assumptions~A.1-A.2 and \cite[Lemma~3.15]{Piet_Jero_HHO_Book_20}, It can be easily shown  that the above seminorms $\|\bullet\|_{a,h}$ and $\|\bullet\|_{1,h}$ are equivalent on $\underline{U}_{h}^k$.

The following boundedness property can be obtained using the definitions of reconstructions $\Gtk,R_T^{k+1}$, trace inequality and \CS inequality, see \cite[Proposition 3.13 \& Lemma 3.15]{Piet_Jero_HHO_Book_20}.
\begin{lemma}[Boundedness]
	For $\underline{u}_h, \underline{v}_h\in\underline{U}_h^k$, there exist a constant $C>0$ independent of mesh parameter $h$ but depending on $\alpha,M,C_{\rm tr}, C_{\rm app}, \varrho,N_{\partial}$ such that
	\begin{align*}
		B_h(\underline{u}_h,\underline{v}_h)\leq C\left(\|\underline{u}_h\|_{a,h}^2+\|u_h\|^2\right)^{1/2} \left(\|\underline{v}_h\|_{a,h}^2+\|v_h\|^2\right)^{1/2}.
	\end{align*}
\end{lemma}

The next lemma is essential for establishing existence and uniqueness of the discrete solution \eqref{HHO_nonself}.
\begin{lemma}[\GardingsType inequality]\label{lem_Garding_ineq}
	For all $\underline{v}_h\in\underline{U}_h^k$, there exists two real numbers $C_1,C_2> 0$ independent of $h$ such that
	\begin{align}\label{garding_ineq}
		B_h(\underline{v}_h,\underline{v}_h)\geq C_1\|\underline{v}_h\|_{a,h}^2- C_2\| v_h\|^2\quad \forall\underline{v}_h\in\underline{U}_h^k.
	\end{align}
\end{lemma}

\begin{proof}

	Following the proof of \cite[Lemma 31.5]{Piet_Jero_HHO_Book_20}, the lower bound for $B_h(\bullet,\bullet)$ can be easily obtained as follows:
	\begin{equation}\label{at_est_coer}
		\left(a_T\Gtk\underline{v}_T, \Gtk\underline{v}_T\right)_T+ s_T(\underline{v}_T,\underline{v}_T)\geq c_3 \|\underline{v}_h\|_{a,h}^2  
	\end{equation}
	for some positive constant $c_3$. Using \CS inequality, we obtain the following estimates for remaining terms
	\begin{align}
		&\sum_{T\in\mathcal{T}_h} \left(\vec{b}{\cdot}\nabla R_T^{k+1}\underline{v}_T, v_T\right)_T\lesssim \|\vec{b}\|_{L^{\infty}(\Omega)} \|\underline{v}_h\|_{a,h} \left(\sum_{T\in\cT_h}\|v_T\|_T^2\right)^{1/2}\leq c_4 \|\underline{v}_h\|_{a,h} \|v_h\|,\label{hho47}\\
		&\sum_{T\in\mathcal{T}_h}(a_0v_T,v_T)_T\leq \|a_0\|_{L^{\infty}(\Omega)} \sum_{T\in\cT_h}\|v_T\|_T^2\leq c_5 \|v_h\|^2,\label{hho48}
	\end{align}
	for some positive constant $c_4$ and $c_5$.
	
	Now, we prove the \GardingsType inequality for $B_h(\underline{v}_h,\underline{v}_h)$. The definition of $B_h(\underline{v}_h,\underline{v}_h)$ and a use of the above three estimates \eqref{at_est_coer}-\eqref{hho48} lead to
	\begin{align*}
		B_h(\underline{v}_h,\underline{v}_h)&\geq c_3\|\underline{v}_h\|_{a,h}^2- c_4\|\underline{v}_h\|_{a,h} \|v_h\|- c_5\|v_h\|^2\\
		&\geq (c_3-c_4\zeta)\|\underline{v}_h\|_{a,h}^2- (\frac{c_4}{\zeta}+c_5)\sum_{T\in\cT_h}\|v_T\|_T^2,
	\end{align*}
	for any positive $\zeta$.
	For sufficiently small choice of $\zeta>0$, the proof of \eqref{garding_ineq} follows.
\end{proof}

Using Lemma~\ref{lem_dis_emb}, we have $\| v_h\|\leq C \|\underline{v}_h\|_{a,h}$. We rewrite the \GardingsType inequality \eqref{garding_ineq} as
\begin{equation}\label{Garding_infsup}
	C_1\|\underline{v}_h\|_{a,h}\leq \sup_{\underline{w}_h\in \underline{U}_{h}^k,\, \|\underline{w}_h\|_{a,h}=1}B_h( \underline{v}_h,\underline{w}_h)+C_2\| v_h\|\fl \underline{v}_h\in \underline{U}_{h}^k.
\end{equation}


\subsection{Existence and uniqueness of discrete solution of nonselfadjoint elliptic problem}
We now show an a priori result of an auxiliary discrete problem which is used to prove the existence and uniqueness of the discrete solution $\underline{u}_h\in \underline{U}_{h,0}^k$ of \eqref{HHO_nonself}.

\begin{lemma}\label{dis_trace7}
	Let $q\in L^2(\Omega)$. Then for sufficiently small $h$, there exists a unique solution $\underline{\phi}_h\in\underline{U}_{h,0}^k$ such that
	\begin{align}
		B_h(\underline{v}_h,\underline{\phi}_h)= (q, v_h)\quad \forall \underline{v}_h\in \underline{U}_{h,0}^k.\label{aux_adj_dis}
	\end{align}
	Moreover, the solution $\underline{\phi}_h$ satisfies
	\begin{align}
		\|\underline{\phi}_h\|_{a,h}\leq C\|q\|.\label{dis_trace8}
	\end{align}
\end{lemma}

\begin{proof}
	We first prove the a priori bound \eqref{dis_trace8}. Since \eqref{aux_adj_dis} leads to a finite dimensional system, the existence and uniqueness of the solution follow from the a priori bound. The \GardingsType inequality \eqref{garding_ineq} with $\underline{v}_h=\underline{\phi}_h$ yields
	\begin{align*}
		C_1\|\underline{\phi}_h\|_{a,h}^2- C_2\|\phi_h\|^2&\leq B_h(\underline{\phi}_h,\underline{\phi}_h).
	\end{align*}
	Using \eqref{aux_adj_dis} and Cauchy--Schwarz inequality, we obtain
	\begin{align*}
		B_h(\underline{\phi}_h,\underline{\phi}_h)= (q,\phi_h)
		&\leq\|q\|_{L^2(\Omega)} \|\phi_h\|_{L^2(\Omega)}\leq (\|q\|^2+\|\phi_h\|^2)/2.    
	\end{align*}
	Combining the above two estimates, we have
	\begin{align}
		\|\underline{\phi}_h\|_{a,h}\leq C_3\|q\|+ C_4\|\phi_h\|.\label{hho52}
	\end{align}
	We apply the Aubin-Nitsche duality argument to estimate $\|\phi_h\|$. For the above $\underline{\phi}_h\in\underline{U}_{h,0}^{k}$, consider the following auxiliary problem:
	\begin{subequations}\label{aux_adj}
		\begin{align}
			-\nabla{\cdot}(a(x)\nabla\psi)+ \vec{b}(x){\cdot}\nabla\psi+a_0(x)\psi&= \phi_h\quad\text{in}~~ \Omega,\\
			\psi&=0\quad\text{on}~~ \partial\Omega.
		\end{align}
	\end{subequations}
	
	The Assumptions A.1-A.2 ensure the existence and uniqueness of the solution $\psi\in H^1_0(\Omega)\cap H^2(\Omega)$ (see \cite[Theorem~3.12]{Ern_Guer_04_FEMs}) and the solution satisfies the following elliptic regularity:
	\begin{align}\label{nonself_apriori}
		\|\psi\|_{H^2(\Omega)}\leq C\|\phi_h\|_{L^2(\Omega)}.
	\end{align}
	Multiply \eqref{aux_adj} by $\phi_h$ and integrate over $\Omega$ to obtain
	\begin{align}\label{IBP_adj_eqn}
		\|\phi_h\|^2= -\integ\nabla{\cdot}(a\nabla\psi)\phi_h\dx+\integ \vec{b}{\cdot}\nabla\psi\phi_h\dx+\integ a_0\psi\phi_h\dx.
	\end{align}
	The integration by parts on the first term of the above equation leads to
	\begin{align}
		&-\integ\nabla{\cdot}(a\nabla\psi)\phi_h\dx=-\sit\nabla{\cdot}(a_T\nabla\psi)\phi_h\dx\notag\\
		&\quad=\sum_{T\in\cT_h}\left(\int_T a_T\nabla \psi{\cdot}\nabla\phi_T\dx+\sum_{F\in\cF_T}\int_F(\phi_F-\phi_T)a_T\nabla\psi{\cdot}\bfnTF\ds\right),\label{diff_integ}
	\end{align}
	where the term related to $\phi_F$ on the skeleton $\cF_h$ is zero owing to the zero boundary condition and \cite[Corollary~1.19]{Piet_Jero_HHO_Book_20}. 
	The local term of the above equation is estimated by the definition of gradient reconstruction \eqref{Grad_recons_proj} as
	\begin{align}
		&\int_T a_T\nabla \psi{\cdot}\nabla\phi_T\dx+\sum_{F\in\cF_T}\int_F(\phi_F-\phi_T)a_T\nabla\psi{\cdot}\bfnTF\ds\notag\\
		&=\int_T a_T\nabla \psi{\cdot}\Gtk\underline{\phi}_T\dx+\sum_{F\in\cF_T}\int_F(\phi_F-\phi_T)\left(a_T\nabla\psi-\pi_T^k(a_T\nabla\psi)\right){\cdot}\bfnTF\ds\notag\\
		&=\int_T a_T\Gtk I_T^k\psi{\cdot}\Gtk\underline{\phi}_T\dx+\int_T a_T(\nabla\psi-\Gtk I_T^k\psi){\cdot}\Gtk\underline{\phi}_T\dx\notag\\
		&\quad+\sum_{F\in\cF_T}\int_F(\phi_F-\phi_T)\left(a_T\nabla\psi-\pi_T^k(a_T\nabla\psi)\right){\cdot}\bfnTF\ds.
	\end{align}
	We rewrite the second term of \eqref{IBP_adj_eqn} as
	\begin{align}
		&\int_\Omega \vec{b}(x){\cdot}\nabla\psi \phi_h \dx=\sit \vec{b}(x){\cdot}\nabla\psi\phi_T\dx\notag\\
		&=\sit \vec{b}(x){\cdot}\nabla(R_T^{k+1} I_T^k \psi)\phi_T\dx+\sit \vec{b}(x){\cdot}\nabla(\psi-R_T^{k+1} I_T^k \psi)\phi_T\dx.\label{hho55}
	\end{align}
	Combining the above equations \eqref{diff_integ}-\eqref{hho55} in \eqref{IBP_adj_eqn} and using the definition of $B_h(\bullet,\bullet)$, we obtain
	\begin{align}
		&\|\phi_h\|^2=B_h(I_h^k\psi,\underline{\phi}_h)-\sum_{T\in\cT_h}s_T(I_T^k\psi,\underline{\phi}_T)+\sumt\Big{(}\int_T a_T(\nabla\psi-\Gtk I_T^k\psi){\cdot}\Gtk\underline{\phi}_T\dx\notag\\
		&\qquad\qquad+\sum_{F\in\cF_T}\int_F(\phi_F-\phi_T)\left(a_T\nabla\psi-\pi_T^k(a_T\nabla\psi)\right){\cdot}\bfnTF\ds\notag\\
		&+\int_T\vec{b}(x){\cdot}\nabla(\psi-R_T^{k+1} I_T^k \psi)\phi_T\dx+\int_T a_0(x)(\psi-\pi_T^k\psi)\phi_h\dx\Big{)}=:T_1+T_2+T_3.\label{aubin_nitche_est}
	\end{align}
	The Cauchy--Schwarz inequality, trace inequality and approximation properties of $R_T^{k+1}I_T^k, \Gtk$ and $\pi_T^k$ yield an estimate for the above last four terms as
	\begin{align}
		T_3\lesssim h\|\psi\|_{H^2(\Omega)}\|\underline{\phi}_h\|_{a,h}.\label{bdd_est}
	\end{align}
	Choose $\underline{v}_h=I_h^k\psi$ in \eqref{aux_adj_dis} and use Cauchy--Schwarz inequality to obtain
	\begin{align}
		B_h(I_h^k\psi,\underline{\phi}_h)= \int_\Omega q\pi_h^k\psi\dx\leq \|q\|\|\pi_h^k\psi\|\lesssim \|q\| \|\psi\|_{H^2(\Omega)}.   
	\end{align}
	Following \cite[equation (46)]{Piet_Ern_Lem_14_arb_local}, we can obtain
	\begin{align}
		\sum_{T\in\cT_h}s_T(I_T^k\psi,\underline{\phi}_T)\lesssim h\|\psi\|_{H^2(\Omega)}\|\underline{\phi}_h\|_{a,h}.\label{stab_est}   
	\end{align}
	Combining last three estimates \eqref{bdd_est}-\eqref{stab_est} in \eqref{aubin_nitche_est} and using the a priori estimate \eqref{nonself_apriori} we obtain
	\begin{align}
		\|\phi_h\|\leq Ch\|\underline{\phi}_h\|_{a,h}+ \|q\|.\label{hho53}
	\end{align}
	Finally, for sufficiently small choice of $h$, use \eqref{hho53} in \eqref{hho52} to obtain
	\begin{align*}
		\|\underline{\phi}_h\|_{a,h}\leq C\|q\|.
	\end{align*}
	This completes the proof.
\end{proof}

\noindent {\bf Existence and uniqueness of the solution of \eqref{HHO_nonself}:}
In Lemma~\ref{dis_trace7}, we proved the existence and uniqueness of the solution of discrete system \eqref{aux_adj_dis} which is the adjoint problem of \eqref{HHO_nonself}. This implies the existence and uniqueness of the solution of discrete system \eqref{HHO_nonself}.


\subsection{Error estimate for nonselfadjoint elliptic problem}
In this subsection, we prove the error estimate for the discrete solution~\eqref{HHO_nonself} using \GardingsType inequality and some auxiliary problems.
\begin{lemma}\label{dis_trace9}
	Let $u\in H_0^1(\Omega)$ and $\underline{u}_h\in \underline{U}_{h,0}^k$ be the solutions to the continuous and discrete problems \eqref{cts_lin_nonselfadj} and \eqref{HHO_nonself}, respectively. Assume $u\in H^{r+2}(\cT_h)$ for some $r\in\{0,1,\ldots,k\}$. For sufficiently small $h$, there exists a real number $C>0$ independent of $h$ such that
	\begin{align}\label{pre_err_est_adj}
		\|I_h^k u- \underline{u}_h\|_{a,h}\leq Ch^{r+1} \|u\|_{H^{r+2}(\cT_h)}.
	\end{align}
\end{lemma}

\begin{proof}
	From the \GardingsType inequality \eqref{garding_ineq}, we have
	\begin{align*}
		C_1\|\underline{v}_h\|_{a,h}^2- C_2\| v_h\|\leq B_h(\underline{v}_h,\underline{v}_h)\quad \forall\underline{v}_h\in\underline{U}_h^k.
	\end{align*}
	Set $\underline{\chi}_h:=I_h^ku-\underline{u}_h$. Choosing $\underline{v}_h=\underline{\chi}_h$ in the above equation, we obtain
	\begin{align}
		C_1\|I_h^ku-\underline{u}_h\|_{a,h}^2- C_2\|\pi_h^ku-u_h\|^2&\leq B_h(I_h^ku-\underline{u}_h,I_h^ku-\underline{u}_h)=B_h(I_h^ku-\underline{u}_h,\underline{\chi}_h).
	\end{align}
	Since $\underline{u}_h\in \underline{U}_{h,0}^k$ satisfies \eqref{HHO_nonself}, we have
	\begin{align}
		C_1\|I_h^ku-\underline{u}_h\|_{a,h}^2- C_2\|\pi_h^ku-u_h\|^2\leq  B_h(I_h^ku,\underline{\chi}_h)-(p,\chi_h).\label{garding_res}
	\end{align}
	Now, $p=-\nabla{\cdot}(a\nabla u)+ \vec{b}{\cdot}\nabla u+a_0u$ of \eqref{cts_lin_nonselfadj} lead to
	\begin{align*}
		(p,\chi_h)&= \sit p\chi_h\dx=\sit(-\nabla{\cdot}(a\nabla u)+ \vec{b}{\cdot}\nabla u+a_0u)\chi_h\dx.
	\end{align*}
	The integration by parts in the above equation and similar approaches from \eqref{IBP_adj_eqn} to \eqref{aubin_nitche_est} yield
	\begin{align*}
		(p,\chi_h)&=B_h(I_h^k u,\underline{\chi}_h)-s_h(I_h^k u,\underline{\chi}_h)+\int_T a_T(\nabla u-\Gtk I_T^k u){\cdot}\Gtk\underline{\chi}_T\dx\notag\\
		&\qquad\qquad+\sum_{F\in\cF_T}\int_F(\chi_F-\chi_T)\left(a_T\nabla u-\pi_T^k(a_T\nabla u)\right){\cdot}\bfnTF\ds\notag\\
		&\qquad+\sit \vec{b}(x){\cdot}\nabla(u-R_T^{k+1} I_T^k u)\chi_T\dx+\integ a_0(x)( u-\pi_T^k u)\chi_h\dx.
	\end{align*}
	The last five terms of the above equation are estimated using similar technique of \eqref{bdd_est} and \eqref{stab_est} and we obtain 
	\begin{align}{\label{res_est_Eh}}
		E_h(p;\underline{\chi}_h):=B_h(I_h^k u,\underline{\chi}_h)-(p,\chi_h)\leq Ch^{r+1} \|u\|_{H^{r+2}(\cT_h)}\|\underline{\chi}_h\|_{a,h}.   
	\end{align}
	This leads to an estimation for \eqref{garding_res} as
	\begin{align}\label{gard_res_est2}
		&C_1\|I_h^k u-\underline{u}_h\|_{a,h}^2- C_2\|\pi_h^ku-u_h\|^2\leq Ch^{r+1} \|u\|_{H^{r+2}(\cT_h)}\|\underline{\chi}_h\|_{a,h}.
	\end{align}
	A use of Lemma~\ref{lem_dis_emb} leads to the above equation as
	\begin{align}
		C_1\|I_h^ku-\underline{u}_h\|_{a,h}\leq Ch^{r+1} \|u\|_{H^{r+2}(\cT_h)}+ C_2\|\chi_h\|.\label{dis_trace4}
	\end{align}
	To estimate $\|\chi_h\|$, we use the results of discrete adjoint problem \eqref{aux_adj_dis} with $q=\chi_h$ and $\underline{v}_h=\underline{\chi}_h$ and this leads to
	\begin{align*}
		\|\chi_h\|^2&= B_h(\underline{\chi}_h,\underline{\phi}_h)=B_h(I_h^ku-\underline{u}_h,\underline{\phi}_h)\notag\\
		&= B_h(I_h^ku,\underline{\phi}_h)-(p,\phi_h)=E_h(p;\underline{\phi}_h).
	\end{align*}
	A use of the estimate for $E_h(p;\underline{\phi}_h)$ in \eqref{res_est_Eh} and a priori result $\|\underline{\phi}_h\|_{a,h}\leq C\|\chi_h\|$ of \eqref{aux_adj_dis}, it holds
	\begin{align}\label{dis_trace2}
		\|\chi_h\|\leq Ch^{r+1} \|u\|_{H^{r+2}(\cT_h)}.
	\end{align}
	Combining \eqref{dis_trace4} and \eqref{dis_trace2}, we obtain the final estimate
	\begin{align*}
		\|I_h^k u-\underline{u}_h\|_{a,h}\leq Ch^{r+1} \|u\|_{H^{r+2}(\cT_h)}.
	\end{align*}
	This completes the proof.
\end{proof}

The triangle inequality, Lemma~\ref{dis_trace9} and Lemma~\ref{lem_apprx_Gtk} lead to the error estimate:
\begin{theorem}[Error estimate]
	Let $u\in H_0^1(\Omega)$ and $\underline{u}_h\in \underline{U}_{h,0}^k$ be the solutions to the continuous and discrete problems \eqref{cts_lin_nonselfadj} and \eqref{HHO_nonself}, respectively. Assume $u\in H^{r+2}(\cT_h)$ for some $r\in\{0,1,\ldots,k\}$. For sufficiently small $h$, it holds
	\begin{align}\label{err_est_nonself_adj}
		\left(\sum_{T\in\cT_h}\|\nabla u-\Gtk\underline{u}_T\|_T^2\right)^{1/2}\leq Ch^{r+1} \|u\|_{H^{r+2}(\cT_h)}
	\end{align}
	for some positive constant $C$ independent of $h$. 
\end{theorem}

\section{Quasilinear elliptic problem}\label{sec:Quasilin}
In this section, we consider HHO approximation for the quasilinear elliptic boundary value problem:
\begin{subequations}\label{quasilin_pde}
	\begin{align}
		-\nabla{\cdot}(a(x,u)\nabla u)&= f(x)\quad\text{in}~~ \Omega,\\
		u(x)&=0\quad\text{on}~~ \partial\Omega,\label{quasi_homo_bdd}
	\end{align}
\end{subequations}
where $\Omega$ is a bounded convex polytopal domain in $\mathbb{R}^d$, $d\in \{2,3\}$ with Lipschitz boundary $\partial\Omega$. 
We make the following assumptions \cite{Doug_Dup_75,Gudi_AKP_07_DG_quasi} for the problem (\ref{quasilin_pde}). 

\medskip

\noindent\textbf{Assumption Q:}
\begin{enumerate}
	\item There exist positive constants $\alpha,M$ such that $0<\alpha\leq a(x, t)\leq M$, $x\in\bar{\Omega},t\in\mathbb{R}$.
	\item The coefficient function $a:\bar{\Omega}\times \mathbb{R}\to\mathbb{R}$ is twice continuously differentiable function on $\bar{\Omega}\times \mathbb{R}$ such that all derivatives of $a(x, t)$ up to and including second order are bounded in $\bar{\Omega}\times \mathbb{R}$, i.e, $a\in C^2_b(\bar{\Omega}\times \mathbb{R})$.
\end{enumerate}

\medskip

Henceforth, we understand $a(x,t)$ by $a(t)$ if there is no confusion. For sufficiently smooth data $f$, the problem \eqref{quasilin_pde}  possesses a unique smooth solution $u$ when the boundary is also sufficiently smooth, see \cite{Dou_Dup_Ser_71_DivExitUniq}. 
Under the assumption $f\in L^{\infty}(\Omega)$, Caloz and Rappaz in \cite[Theorem 5.1]{Caloz_Rappaz_97_Nonlin_Bifur} have proved that the solution of \eqref{quasilin_pde} has unique solution in $W^{2,p}(\Omega),\, 1\leq p <\infty$ on a regular domain $\Omega$. The existence and uniqueness of the weak solution $u$ in $H^1_0(\Omega)$ is established in \cite{Ivan_Michal_94_WeakSoln_Quasilin} on bounded domain with Lipschitz boundary.

Finite element approximations for the above problem \eqref{quasilin_pde} is proposed and analyzed under various assumptions on the coefficient and data, and assumption on regularity for the solution, see \cite{Doug_Dup_75,Lip_Michal_96,Milner_85,Gudi_AKP_07_DG_quasi,Chun_Vict_11_TwoGridDG_Quasi,Chun_Vict_13_Apost_DG_Quasi,Chun_Vict_09_APost_FV_Quasi,Chun_Vict_TwoGrid_FV}. Following \cite{Canggiani_20_VEM_Quasilin,Gudi_NN_AKP_08_HpDG_Quasi,Chun_Vict_11_TwoGridDG_Quasi,Chun_Vict_13_Apost_DG_Quasi}, we make the following assumption on the solution of \eqref{quasilin_pde} as:

\medskip

\noindent\textbf{Assumption Q}.3. Assume the solution $u \in H^1_0(\Omega)$ of \eqref{quasilin_pde} belongs to $H^2(\Omega)\cap W^{1,\infty}(\Omega)$ for $d=2$ and belongs to $H^3(\Omega)$ for $d=3$.

\medskip

\begin{remark}
	For our subsequent error analysis, Assumption Q.3 can be relaxed to $u\in H^{3/2+\epsilon}(\Omega)\cap W^{1,\infty}(\Omega),\,\epsilon>0$ for $d=2$, and to $u\in H^1_0(\Omega)\cap H^{5/2+\epsilon}(\Omega),\,\epsilon>0$ for $d=3$. However, these require the approximation properties of \eqref{proj_est} and \eqref{eqn_recons_est} related to the projections $\pi_T^k$ and $\pi_F^k$ on fractional order Sobolev spaces, see \cite[Remark~1.49]{Piet_Jero_HHO_Book_20}. For simplicity of presentation, we kept our assumptions on integral Sobolev spaces.
\end{remark}

For the HHO approximation of the solution of \eqref{quasilin_pde}, we consider a linearized problem: find $\psi\in H^1_0(\Omega)$ such that
\begin{subequations}\label{lin_nonself}
	\begin{align}
		-\nabla{\cdot}(a(u)\nabla\psi+a_u(u)\nabla u\,\psi)&= \phi\quad\text{in}~~ \Omega,\\
		\psi&=0\quad\text{on}~~ \partial\Omega,
	\end{align}
\end{subequations}
for some load function $\phi\in\lt$.
Since $a\in C^2_b(\mathbb{R})$ and $u\in W^{1,\infty}(\Omega)$, the existence and uniqueness of the solution of above linear elliptic problem \eqref{lin_nonself} follow. Moreover, the solution satisfies the elliptic regularity $\|\psi\|_{H^2(\Omega)}\leq C\|\phi\|$, see \cite[Theorem~3.12]{Ern_Guer_04_FEMs}. 
In the following sections, we consider a HHO approximation of the above linearized problem~\eqref{lin_nonself} and analyze the existence and uniqueness of the HHO approximation for the quasilinear problem \eqref{quasilin_pde}.

\subsection{HHO approximation for quasilinear elliptic problem}
For $\underline{w}_h, \underline{u}_h$ and $\underline{v}_h$ in $\underline{U}_{h}^k$, define
\begin{align}\label{defn_Bh}
	\mathcal{N}_h(\underline{w}_h;\underline{u}_h, \underline{v}_h)&:=\sum_{T\in\mathcal{T}_h} \int_T
	a(R_T^{k+1}\underline{w}_T) \Gtk\underline{u}_T{{\cdot}} \Gtk\underline{v}_T\dx+ s_h(\underline{u}_h,\underline{v}_h),
\end{align}
where the above stabilization term $s_h(\underline{u}_h,\underline{v}_h):=\sum_{T\in\mathcal{T}_h}s_T(\underline{u}_T,\underline{v}_T)$ with the local stabilization 
\begin{align*}
	s_T(\underline{u}_T,\underline{v}_T):= \sum_{F\in\mathcal{F}_T}\frac{a_{TF}^{\infty}}{h_F} \left(S_F^k\underline{u}_T,S_F^k\underline{v}_T\right)_F \text{ where } a_{TF}^{\infty}:=\|a_T\|_{L^{\infty}(F)}.
\end{align*}
Discrete hybrid-high order approximation of \eqref{quasilin_pde} seeks $\underline{u}_h\in \underline{U}_{h,0}^k$ such that
\begin{align}\label{hho_dis_quasi}
	\mathcal{N}_h(\underline{u}_h;\underline{u}_h, \underline{v}_h)=l(\underline{v}_h)\quad \forall\underline{v}_h\in \underline{U}_{h,0}^k,
\end{align}
where the linear form reads
$l(\underline{v}_h):=\sum_{T\in\mathcal{T}_h} \int_T fv_T\dx$.

We propose a hybrid high-order discretization of linearized problem \eqref{lin_nonself}.
The discrete linearized problem seeks $\underline{\psi}_h\in \underline{U}_{h,0}^k$ such that
\begin{align}\label{dis_lin_quasi}
	\mathcal{N}_h^{\text{lin}}(u;\underline{\psi}_h, \underline{v}_h)=(\phi,v_h)\quad \forall\underline{v}_h\in \underline{U}_{h,0}^k,
\end{align}
where, we considered linearization around the solution $u$ of \eqref{quasilin_pde} and for $\underline{\psi}_h, \underline{v}_h\in \underline{U}_h^k$,
\begin{align}
	\mathcal{N}_h^{\text{lin}}(u; \underline{\psi}_h, \underline{v}_h)&:= \sum_{T\in\mathcal{T}_h}\int_T \left(a(u) \Gtk\underline{\psi}_T{\cdot} \Gtk\underline{v}_T +a_u(u)R_T^{k+1}\underline{\psi}_T\nabla u{\cdot} \Gtk\underline{v}_T \right)\dx\notag\\
	&\qquad+s_h( \underline{\psi}_h,\underline{v}_h).\label{defn_Bh_lin}
\end{align}
Moreover, define a fully discrete version of linearization term as: for $\underline{w}_h,\underline{\psi}_h, \underline{v}_h\in \underline{U}_h^k$,
\begin{align}
	\widetilde{\mathcal{N}}_h^{\text{lin}}(\underline{w}_h; \underline{\psi}_h, \underline{v}_h)&:= \sum_{T\in\mathcal{T}_h}\int_T a(R_T^{k+1}\underline{w}_T) \Gtk\underline{\psi}_T{\cdot} \Gtk\underline{v}_T\dx\notag\\ &\quad+\sit a_u(R_T^{k+1}\underline{w}_T)R_T^{k+1}\underline{\psi}_T\Gtk\underline{w}_T{\cdot} \Gtk\underline{v}_T\dx
	+s_h( \underline{\psi}_h,\underline{v}_h).\label{defn_Bh_lin_dis}
\end{align}
For $\underline{v}_h\in \underline{U}_h^k$, the definition of reconstruction operator yields $\int_T (R_T^{k+1}\underline{v}_T -v_T)\dx=0$ for each $T\in\cT_h$.  This leads to the following estimate (see \cite[Corollary 5.10]{Piet_Dron_16_Leray}): for $p\leq q\leq p^*$,
\begin{align}
	\|R_h^{k+1}\underline{v}_h-v_h\|_{L^q(\Omega)}\leq C h^{1+\frac{d}{q}-\frac{d}{p}}\|  \underline{v}_h\|_{1,p,h}.\label{recons_est_Lp} 
\end{align}
In particular, 
\begin{equation}\label{L6_bdd_H1}
	\|R_h^{k+1}\underline{v}_h\|_{L^6(\Omega)}\leq C\|\underline{v}_h\|_{1,h}.
\end{equation}
Using similar arguments as in the proof of Lemma~\ref{lem_Garding_ineq} and the above estimate, we obtain a \GardingsType inequality for $\mathcal{N}_h^{\text{lin}}(u; \bullet, \bullet)$ as
\begin{align}
	\mathcal{N}_h^{\text{lin}}(u; \underline{v}_h, \underline{v}_h)\geq C_1\|\underline{v}_h\|_{1,h}^2-C_2\|v_h\|_{L^2(\Omega)}^2   
\end{align}
for some positive constants $C_1$ and $C_2$ independent of $h$ but depending on $u$ and $a(u)$. Then, the existence and uniqueness of the solution $\underline{\psi}_h\in \underline{U}_{h,0}^k$ of \eqref{dis_lin_quasi} follow from the existence and uniqueness of the solution of \eqref{HHO_nonself}.

We make the following assumption \cite{Gudi_AKP_07_DG_quasi} throughout the section.

\noindent\textbf{Assumption Q}.4 (Quasi-uniformity). We assume the admissible mesh sequence $(\cT_h)_{h>0}$ to be quasi-uniform, i.e., there exists a constant $C_Q$ independent of $h$ such that 
\begin{equation*}
	\max_{T\in\cT_h} h_T\leq C_Q\min_{T\in\cT_h}h_T.
\end{equation*}

\subsection{Fixed point formulation and contraction result}
In this section, we prove the existence, local uniqueness and error estimates for the solution $\underline{u}_h\in \underline{U}_{h,0}^k$ of the above problem \eqref{hho_dis_quasi} using fixed point arguments. Following the idea of \cite{GM_NN_16_VKE_NCFEM,CC_GM_NN_19_VKE_DG,Gudi_AKP_07_DG_quasi}, we define a nonlinear map $\mu: \underline{U}_{h,0}^k\to \underline{U}_{h,0}^k$ which satisfies, for all $\underline{v}_h\in \underline{U}_{h,0}^k$
\begin{align}
	\mathcal{N}_h^{\text{lin}}(u;I_h^ku-\mu(\underline{\theta}_h), \underline{v}_h)= \mathcal{N}_h^{\text{lin}}(u;I_h^ku-\underline{\theta}_h, \underline{v}_h)+ \mathcal{N}_h(\underline{\theta}_h;\underline{\theta}_h, \underline{v}_h)-l(\underline{v}_h).\label{nonlinear_map_mu}
\end{align}
The map $\mu$ is well-defined as \eqref{dis_lin_quasi} is well-posed. It can be observed that any fixed point $\underline{\xi}_h$ (say) of $\mu$ satisfies $\mathcal{N}_h(\underline{\xi}_h;\underline{\xi}_h, \underline{v}_h)=l(\underline{v}_h)$ for all $\underline{v}_h\in \underline{U}_{h,0}^k$. That is, $\underline{\xi}_h$ is a solution of \eqref{hho_dis_quasi}. Define a ball of radius $R$ with center at $I_h^k u$ as
\begin{equation*}
	D(I_h^ku;R):=\left\{\underline{\theta}_h\in\underline{U}_{h,0}^k \text{ such that } \|I_h^k u-\underline{\theta}_h\|_{1,h}\leq R\right\}.
\end{equation*}

The following results are obtained using generalized \Holders inequality, Lemma~\ref{lem_Sob_Inv_Ineq}, equation~\eqref{L6_bdd_H1}, discrete Sobolev embedding \eqref{dis_emb} and similar arguments of \cite[Equation~4.20]{Gudi_AKP_07_DG_quasi} as:
\begin{lemma}\label{lem_imp_ineq}
	For $\underline{\xi}_h,\underline{\chi}_h\in \underline{U}_{h}^k$ and $\underline{v}_h\in \underline{U}_{h,0}^k$, the following bounds hold true
	\begin{enumerate}[(a)]
		\item $\displaystyle\sum_{T\in\cT_h} \int_T |R_T^{k+1}\underline{\xi}_T\Gtk\underline{\chi}_T{\cdot}\Gtk\underline{v}_T|\dx\\ 
		\hspace*{5mm}\leq C\big{(}\max_{T\in\cT_h}h_T^{-d/6}\big{)} \|R_h^{k+1}\underline{\xi}_h\|_{L^6(\Omega)}\|\underline{\chi}_h\|_{1,h}\|\underline{v}_h\|_{1,h}\\
		\hspace*{5mm}\leq C\big{(}\max_{T\in\cT_h}h_T^{-d/6}\big{)} \|\underline{\xi}_h\|_{1,h}\|\underline{\chi}_h\|_{1,h}\|\underline{v}_h\|_{1,h}.$
		\item $\displaystyle\sum_{T\in\cT_h} \int_T |(R_T^{k+1}\underline{\xi}_T)^2\Gtk\underline{\chi}_T{\cdot}\Gtk\underline{v}_T|\dx\\
		\hspace*{5mm}\leq C\big{(}\max_{T\in\cT_h}h_T^{-d/3}\big{)} \|R_h^{k+1}\underline{\xi}_h\|_{L^6(\Omega)}^2\|\underline{\chi}_h\|_{1,h}\|\underline{v}_h\|_{1,h}\\
		\hspace*{5mm}\leq C\big{(}\max_{T\in\cT_h}h_T^{-d/3}\big{)} \|\underline{\xi}_h\|_{1,h}^2\|\underline{\chi}_h\|_{1,h}\|\underline{v}_h\|_{1,h}.$
	\end{enumerate}
\end{lemma}

Now we prove some auxiliary results which will be used in the fixed point theorem.
\begin{lemma}\label{lem:diff_lin}
	Let $u\in H^{r+2}(\cT_h)$ for $r\in\{0,1,\ldots,k\}$ and $\underline{\psi}_h, \underline{v}_h\in \underline{U}_{h}^k$, it holds
	\begin{align}
		|\mathcal{N}_h^{\text{lin}}(u; \underline{\psi}_h, \underline{v}_h)-\widetilde{\mathcal{N}}_h^{\text{lin}}(I_h^k u; \underline{\psi}_h, \underline{v}_h)|\leq Ch^{r+1-d/6}\|u\|_{H^{r+2}(\cT_h)}\|\underline{\psi}_h\|_{1,h}\|\underline{v}_h\|_{1,h}.\label{eq_est_diff_lin}
	\end{align}
\end{lemma}
\begin{proof}
	From the definitions of $\mathcal{N}_h^{\text{lin}}$ and $\widetilde{\mathcal{N}}_h^{\text{lin}}$ in \eqref{defn_Bh_lin}-\eqref{defn_Bh_lin_dis}, we have
	\begin{align}
		&\mathcal{N}_h^{\text{lin}}(u; \underline{\psi}_h, \underline{v}_h)-\widetilde{\mathcal{N}}_h^{\text{lin}}(I_h^k u; \underline{\psi}_h, \underline{v}_h)\notag\\
		&= \sum_{T\in\mathcal{T}_h}\int_T (a(u)- a(R_T^{k+1}I_T^k u))\Gtk\underline{\psi}_T{\cdot} \Gtk\underline{v}_T\dx\notag\\
		&\quad+\sum_{T\in\mathcal{T}_h}\int_T \left(a_u(u)\nabla u-a_u(R_T^{k+1}I_T^k u)\Gtk I_T^k u\right)R_T^{k+1}\underline{\psi}_T{\cdot} \Gtk\underline{v}_T\dx.\label{eq_diff_lin}
	\end{align}
	The Taylor series expansion 
	\begin{align}\label{Taylor_Expan_au}
		a(u)=a(w)+ \tilde{a}_u(u)(u-w),    
	\end{align}
	where $\tilde{a}_{u}(u)=\int_0^1a_{u}(u+t(w-u))\dt$, \Holders inequality, inverse inequality and Lemma~\ref{lem_apprx_recons} lead to an estimate for the first term of \eqref{eq_diff_lin} as
	\begin{align}
		&\sum_{T\in\mathcal{T}_h}\int_T (a(u)- a(R_T^{k+1}I_T^k u))\Gtk\underline{\psi}_T{\cdot} \Gtk\underline{v}_T\dx\notag\\
		&\leq C \|\tilde{a}_u\|_{L^\infty(\Omega)}\sum_{T\in\mathcal{T}_h}\|u-R_h^{k+1}I_h^k u\|_{L^2(T)}  \|\Ghk\underline{\psi}_h\|_{L^4(T)}\|\Ghk\underline{v}_h\|_{L^4(T)}\notag\\
		&\leq C\|\tilde{a}_u\|_{L^\infty(\Omega)}h^{r+2-d/2}\|u\|_{H^{r+2}(\cT_h)}\|\underline{\psi}_h\|_{1,h}\|\underline{v}_h\|_{1,h}.
	\end{align}
	The second term of \eqref{eq_diff_lin} is estimated using Taylor series expansion, \Holders inequality, inverse inequality, quasi-uniformity of meshes,  Lemma~\ref{lem_dis_emb}-\ref{lem_apprx_recons} and Lemma~\ref{lem_imp_ineq}(a) as
	\begin{align}
		&\sum_{T\in\mathcal{T}_h}\int_T \left(a_u(u)\nabla u-a_u(R_T^{k+1}I_T^k u)\Gtk I_T^k u\right)R_T^{k+1}\psi_T{\cdot} \Gtk\underline{v}_T\dx\notag\\
		&=\sum_{T\in\mathcal{T}_h}\int_T \left(a_u(u)-a_u(R_T^{k+1}I_T^k u)\right)(\nabla u) R_T^{k+1}\psi_T{\cdot} \Gtk\underline{v}_T\dx\notag\\
		&\quad+\sum_{T\in\mathcal{T}_h}\int_T a_u(R_T^{k+1}I_T^k u)\left(\nabla u-\Gtk I_T^k u\right)R_T^{k+1}\psi_T{\cdot} \Gtk\underline{v}_T\dx\notag\\
		&\leq C\|\tilde{a}_{uu}\|_{L^\infty(\Omega)}h^{r+2-d/6}\|u\|_{H^{r+2}(\cT_h)}\|\nabla u\|_{L^\infty(\Omega)}\|\underline{\psi}_h\|_{1,h}\|\underline{v}_h\|_{1,h}\notag\\
		&\quad+C\|\tilde{a}_{u}\|_{L^\infty(\Omega)}h^{r+1-d/6}\|u\|_{H^{r+2}(\cT_h)}\|\underline{\psi}_h\|_{1,h}\|\underline{v}_h\|_{1,h}.
	\end{align}
	Then, the proof of \eqref{eq_est_diff_lin} follows from the above two estimations.
\end{proof}

\begin{theorem}[Fixed point result]\label{thm_fixed_point}
	Let $u\in H^1_0(\Omega)$ be a solution for \eqref{quasilin_pde}. Assume $u\in H^{r+2}(\cT_h)$ for some $r\in\{d-2,\ldots,k\}$. For sufficiently small mesh parameter $h$, there exists $R(h)$ such that the nonlinear map $\mu: \underline{U}_{h,0}^k\to \underline{U}_{h,0}^k$ defined in \eqref{nonlinear_map_mu} maps from the ball $D(I_h^ku;R(h))$ to itself with radius $R(h):=\tilde{C}h^{r+1}$ for some constant $\tilde{C}$ independent of the mesh parameter. Moreover, $\mu$ has a fixed point in $D(I_h^ku;R(h))$.
\end{theorem}

\begin{proof}
	Since $\mathcal{N}_h^{\text{lin}}(u;\bullet,\bullet)$ is associated with the linearized problem, it satisfies the \GardingsType inequality \eqref{Garding_infsup}:
	\begin{align}
		C_1\|\underline{\theta}_h\|_{1,h}\leq \mathcal{N}_h^{\text{lin}}(u;\underline{\theta}_h, \underline{v}_h)+ C_2\|\theta_h\|\fl \underline{\theta}_h\in\underline{U}_{h,0}^k\label{Garding_infsup_2}
	\end{align}
	for some $\underline{v}_h\in\underline{U}_{h,0}^k$ with $\|\underline{v}_h\|_{1,h}=1$. Replacing $\underline{\theta}_h$ with $I_h^ku-\mu(\underline{\theta}_h)$ (we understand $\|\underline{\psi}_h\|_{L^2}$ by $\|\psi_h\|_{L^2(\Omega)}$) and using the definition of $\mu$ of \eqref{nonlinear_map_mu}, we obtain
	\begin{align}
		&C_1\|I_h^ku-\mu(\underline{\theta}_h)\|_{1,h}\leq \mathcal{N}_h^{\text{lin}}(u; I_h^ku-\mu(\underline{\theta}_h), \underline{v}_h)+ C_2\|I_h^ku-\mu(\underline{\theta}_h)\|_{L^2}\notag\\
		& =\mathcal{N}_h^{\text{lin}}(u;I_h^ku-\underline{\theta}_h, \underline{v}_h)+ \mathcal{N}_h(\underline{\theta}_h;\underline{\theta}_h, \underline{v}_h)-l(\underline{v}_h)+ C_2\|I_h^ku-\mu(\underline{\theta}_h)\|_{L^2}\notag\\
		& =\left(\mathcal{N}_h^{\text{lin}}(u;I_h^ku-\underline{\theta}_h, \underline{v}_h)-\widetilde{\mathcal{N}}_h^{\text{lin}}(I_h^ku;I_h^ku-\underline{\theta}_h, \underline{v}_h)\right)
		\notag\\
		&\quad+\left(\widetilde{\mathcal{N}}_h^{\text{lin}}(I_h^ku;I_h^ku-\underline{\theta}_h, \underline{v}_h)+ \mathcal{N}_h(\underline{\theta}_h;\underline{\theta}_h, \underline{v}_h)-l(\underline{v}_h)\right)+ C_2\|I_h^ku-\mu(\underline{\theta}_h)\|_{L^2}.\label{garding_est_1}
	\end{align}
	The definition of $l(\underline{v}_h)$ with \eqref{quasilin_pde} and integration by parts lead to 
	\begin{align}
		l(\underline{v}_h)&=\integ fv_h\dx=-\integ\nabla{\cdot}(a(u)\nabla u) v_h\dx=-\sit\nabla{\cdot}(a(u)\nabla u) v_h\dx\notag\\
		&\quad=\sum_{T\in\cT_h}\left(\int_T a(u)\nabla  u{\cdot}\nabla v_T\dx+\sum_{F\in\cF_T}\int_F( v_F- v_T)a(u)\nabla u{\cdot}\bfnTF\ds\right),
	\end{align}
	where the additional term related to $v_F$ on the skeleton $\cF_h$ is zero, owing to the zero boundary condition and \cite[Corollary~1.19]{Piet_Jero_HHO_Book_20}.
	Now we rewrite the above terms by adding and subtracting several terms and use the definition \eqref{Grad_recons_proj} of Gradient reconstruction as follows
	\begin{align}
		&l(\underline{v}_h)=\sum_{T\in\cT_h}\left(\int_T a(u)\nabla  u{\cdot}\Gtk \underline{v}_T\dx+\sum_{F\in\cF_T}\int_F( v_F- v_T)(a(u)\nabla u-\pi_T^k(a(u)\nabla u)){\cdot}\bfnTF\ds\right)\notag\\
		&=\sum_{T\in\cT_h}\int_T a(u)\Gtk I_T^k u{\cdot}\Gtk \underline{v}_T\dx+\sum_{T\in\cT_h}\int_T a(u)(\nabla  u-\Gtk I_T^k u){\cdot}\Gtk \underline{v}_T\dx\notag\\
		&\qquad+\sum_{T\in\cT_h}\sum_{F\in\cF_T}\int_F( v_F- v_T)(a(u)\nabla u-\pi_T^k(a(u)\nabla u)){\cdot}\bfnTF\ds\notag\\
		&=\sum_{T\in\cT_h}\int_T a(R_T^{k+1}I_T^k u)\Gtk I_T^k u{\cdot}\Gtk \underline{v}_T\dx+\sum_{T\in\cT_h}\int_T \left(a(u)-a(R_T^{k+1}I_T^k u)\right)\Gtk I_T^k u{\cdot}\Gtk \underline{v}_T\dx\notag\\
		&\qquad+\sum_{T\in\cT_h}\int_T a(u)(\nabla  u-\Gtk I_T^k u){\cdot}\Gtk\underline{v}_T\dx\notag\\
		&\quad\qquad+\sum_{T\in\cT_h}\sum_{F\in\cF_T}\int_F( v_F- v_T)(a(u)\nabla u-\pi_T^k(a(u)\nabla u)){\cdot}\bfnTF\ds.\label{nonlin_res}
	\end{align}
	The second term of the above equation~\eqref{nonlin_res} is estimated using Taylor series expansion \eqref{Taylor_Expan_au}, generalized \Holders inequality and Lemma~\ref{lem_apprx_recons} as
	\begin{align}
		&\sum_{T\in\cT_h}\int_T \left(a(u)-a(R_T^{k+1}I_T^k u)\right)\Gtk I_T^k u{\cdot}\Gtk \underline{v}_T\dx\leq Ch^{r+2}\|u\|_{H^{r+2}(\cT_h)}\|u\|_{W^{1,\infty}(\Omega)}\|\underline{v}_h\|_{1,h}.\label{au_proj_res}
	\end{align}
	The third term of \eqref{nonlin_res} is estimated using \CS inequality and Lemma~\ref{lem_apprx_Gtk} as
	\begin{align}
		\sum_{T\in\cT_h}\int_T a(u)(\nabla  u-\Gtk I_T^k u){\cdot}\Gtk\underline{v}_T\dx\leq Ch^{r+1}\|u\|_{H^{r+2}(\cT_h)}\|\underline{v}_h\|_{1,h}.
	\end{align}
	The last term of \eqref{nonlin_res} is estimated by the Cauchy--Schwarz inequality, Lemma~\ref{proj_est}, trace inequality and Sobolev embedding $H^1(\Omega)\hookrightarrow L^4(\Omega)$ as
	\begin{align}
		&\sum_{T\in\cT_h}\sum_{F\in\cF_T}\int_F( v_F- v_T)\left(a(u)\nabla u-\pi_T^k (a(u)\nabla  u)\right){\cdot}\bfnTF\ds\notag\\
		&\quad\leq Ch^{r+1}\|a(u)\nabla u\|_{H^{r+1}(\cT_h)}\|\underline{v}_h\|_{1,h}\notag\\
		&\quad\leq C h^{r+1} (\|a'(u)\|_{L^\infty}\|u\|_{H^{r+2}(\cT_h)}+\|a(u)\|_{L^\infty})\|u\|_{H^{r+2}(\cT_h)}\|\underline{v}_h\|_{1,h}. \label{au_proj}
	\end{align}

	For $\underline{\xi}_h, \underline{v}_h\in\underline{U}_{h}^k$, define
	\begin{equation}\label{defn_Fh}
		\langle \mathcal{F}_h(\underline{\xi}_h),\underline{v}_h\rangle:= \sum_{T\in\mathcal{T}_h}\int_T a(R_T^{k+1}\underline{\xi}_T) \Gtk\underline{\xi}_T{\cdot} \Gtk\underline{v}_T\dx.    
	\end{equation}
	The definitions of $\mathcal{N}_h(\bullet;\bullet,\bullet)$ and $\mathcal{F}_h$, and previous estimates \eqref{au_proj}-\eqref{au_proj_res} lead to
	\begin{align}
		&\mathcal{N}_h(\underline{\theta}_h;\underline{\theta}_h, \underline{v}_h)-l(\underline{v}_h)\leq \langle \mathcal{F}_h(\underline{\theta}_h),\underline{v}_h\rangle -\langle \mathcal{F}_h(I_h^k u),\underline{v}_h\rangle +s_h(\underline{\theta}_h,\underline{v}_h)+Ch^{r+1}\|\underline{v}_h\|_{1,h}.\label{res_est_1}   
	\end{align}
	The Taylor series expansion yields
	\begin{align}
		a(w)= a(u)+ a_u(u)(w-u)+ \tilde{a}_{uu}(w)(w-u)^2,\label{Taylor_Expan}
	\end{align}
	where $\tilde{a}_{uu}(w)=\int_0^1(1-t)a_{uu}(w+t(w-u))\dt$. Since $a_u\in C^1_b(\bar{\Omega}\times \bR)$ and $a_{uu}\in C^0_b(\bar{\Omega}\times \bR)$, we have $\tilde{a}_u\in L^{\infty}(\Omega\times \bR)$ and $\tilde{a}_{uu}\in L^{\infty}(\Omega\times \bR)$, see \cite[equation (4.8)]{Gudi_AKP_07_DG_quasi}. We set
	\begin{equation}
		C_a:=\max\left\{\|\tilde{a}_u\|_{L^{\infty}(\Omega\times \bR)},\|\tilde{a}_{uu}\|_{L^{\infty}(\Omega\times \bR)} \right\}.  \end{equation}
	For $\underline{\xi}_h\in\underline{U}_{h}^k$ and $\underline{\chi}_h,\underline{v}_h \in\underline{U}_{h,0}^k$, expanding $\mathcal{F}_h(\underline{\xi}_h+\underline{\chi}_h)$ from the above definition \eqref{defn_Fh} and using \eqref{Taylor_Expan}, we obtain
	\begin{align}
		&\langle \mathcal{F}_h(\underline{\xi}_h+\underline{\chi}_h),\underline{v}_h\rangle= \sit a(R_T^{k+1}(\underline{\xi}_T+\underline{\chi}_T)) \Gtk(\underline{\xi}_T+\underline{\chi}_T){\cdot} \Gtk\underline{v}_T\dx\notag\\
		&= \sit \left(a(R_T^{k+1}\underline{\xi}_T)+ a_u(R_T^{k+1}\underline{\xi}_T)R_T^{k+1}\underline{\chi}_T\right)\Gtk(\underline{\xi}_T+\underline{\chi}_T){\cdot} \Gtk\underline{v}_T\dx\notag\\
		&\qquad+ \sit\tilde{a}_{uu}(R_T^{k+1}\underline{\xi}_T)(R_T^{k+1}\underline{\chi}_T)^2\Gtk(\underline{\xi}_T+\underline{\chi}_T){\cdot} \Gtk\underline{v}_T\dx\notag\\
		&= \langle \mathcal{F}_h(\underline{\xi}_h),\underline{v}_h\rangle+\widetilde{\mathcal{N}}_h^{\text{lin}}(\underline{\xi}_h;\underline{\chi}_h, \underline{v}_h)-s_h(\underline{\chi}_h,\underline{v}_h)\notag\\
		&\quad+ \sit\tilde{a}_{uu}(R_T^{k+1}\underline{\xi}_T)(R_T^{k+1}\underline{\chi}_T)^2\Gtk\underline{\xi}_T{\cdot} \Gtk\underline{v}_T\dx\notag\\
		&\qquad+\sit a_u(R_T^{k+1}\underline{\xi}_T)R_T^{k+1}\underline{\chi}_T\Gtk\underline{\chi}_T{\cdot}\Gtk\underline{v}_T\dx\notag\\
		&\quad\qquad+ \sit \tilde{a}_{uu}(R_T^{k+1}\underline{\xi}_T)(R_T^{k+1}\underline{\chi}_T)^2\Gtk\underline{\chi}_T{\cdot} \Gtk\underline{v}_T\dx.\label{Gh_expan}  
	\end{align}
	The fifth term of the above equation \eqref{Gh_expan} is estimated using Lemma~\ref{lem_imp_ineq}(a) as
	\begin{align}
		&\sit a_u(R_T^{k+1}\underline{\xi}_T)R_T^{k+1}\underline{\chi}_T\Gtk\underline{\chi}_T{\cdot}\Gtk\underline{v}_T\dx\notag\\
		&\leq C_a C \big{(}\max_{T\in\cT_h}h_T^{-d/6}\big{)}\|\underline{\chi}_h\|_{1,h}^2\|\underline{v}_h\|_{1,h}.
	\end{align}
	The fourth term of \eqref{Gh_expan} is estimated using Lemma~\ref{lem_imp_ineq}(b) as
	\begin{align}
		&\sit\tilde{a}_{uu}(R_T^{k+1}\underline{\xi}_T)(R_T^{k+1}\underline{\chi}_T)^2\Gtk\underline{\xi}_T{\cdot} \Gtk\underline{v}_T\dx\notag\\
		&\leq C_a C\big{(}\max_{T\in\cT_h}h_T^{-d/3}\big{)}\|\underline{\chi}_h\|_{1,h}^2\|\underline{\xi}_h\|_{1,h}\|\underline{v}_h\|_{1,h}\label{est_typical_2}
	\end{align}
	and the sixth term of \eqref{Gh_expan} is estimated by \eqref{est_typical_2} with $\underline{\xi}_h=\underline{\chi}_h$. 
	Combining \eqref{Gh_expan}-\eqref{est_typical_2} with $\underline{\chi}_h:=\underline{\theta}_h-I_h^k u$ and $\underline{\xi}_h=I_h^k u$, we obtain
	\begin{align}
		&\langle \mathcal{F}_h(\underline{\theta}_h),\underline{v}_h\rangle -\langle \mathcal{F}_h(I_h^k u),\underline{v}_h\rangle +s_h(\underline{\theta}_h,\underline{v}_h)\notag\\
		&\leq \widetilde{\mathcal{N}}_h^{\text{lin}}(I_h^k u;\underline{\theta}_h-I_h^k u, \underline{v}_h)+s_h(I_h^k u,\underline{v}_h)\notag\\
		&+C_a C\max_{T\in\cT_h}h_T^{-d/3}\left(\|\underline{\theta}_h-I_h^k u\|_{1,h}^2+\|\underline{\theta}_h-I_h^k u\|_{1,h}^2\|I_h^k u\|_{1,h}+\|\underline{\theta}_h-I_h^k u\|_{1,h}^3\right)\|\underline{v}_h\|_{1,h}\notag\\
		&\leq \widetilde{\mathcal{N}}_h^{\text{lin}}(I_h^k u;\underline{\theta}_h-I_h^k u, \underline{v}_h)+C_a C\max_{T\in\cT_h}h_T^{-d/3}\left(\|\underline{\theta}_h-I_h^k u\|_{1,h}^2+\|\underline{\theta}_h-I_h^k u\|_{1,h}^3\right)\|\underline{v}_h\|_{1,h}\notag\\
		&\qquad+Ch^{r+1}\|\underline{v}_h\|_{1,h},\label{Gh_err_stab}
	\end{align}
	where we used the estimation for the consistency term $s_h(I_h^k u,\underline{v}_h)$, see \cite[equation (46)]{Piet_Ern_Lem_14_arb_local}.
	Combining \eqref{res_est_1} and \eqref{Gh_err_stab}, we obtain
	\begin{align}
		&\widetilde{\mathcal{N}}_h^{\text{lin}}(I_h^ku;I_h^ku-\underline{\theta}_h, \underline{v}_h)+ \mathcal{N}_h(\underline{\theta}_h;\underline{\theta}_h, \underline{v}_h)-l(\underline{v}_h)\notag\\
		&\leq C_a C\big{(}\max_{T\in\cT_h}h_T^{-d/3}\big{)}\left(\|\underline{\theta}_h-I_h^k u\|_{1,h}^2+\|\underline{\theta}_h-I_h^k u\|_{1,h}^3\right)\|\underline{v}_h\|_{1,h}+Ch^{r+1}\|\underline{v}_h\|_{1,h}.\label{com_res_est}
	\end{align}
	This implies from \eqref{garding_est_1} with Lemma~\ref{lem:diff_lin} and  $\|\underline{v}_h\|_{1,h}=1$ that
	\begin{align}
		C_1\|I_h^ku-\mu(\underline{\theta}_h)\|_{1,h}&\leq C\Big{(}C_a \big{(}\max_{T\in\cT_h}h_T^{-d/3}\big{)}\left(\|\underline{\theta}_h-I_h^k u\|_{1,h}^2+\|\underline{\theta}_h-I_h^k u\|_{1,h}^3\right)\notag\\
		&\quad+h^{r+1}+h^{r+1-d/6}\|\underline{\theta}_h-I_h^k u\|_{1,h}+ \|I_h^ku-\mu(\underline{\theta}_h)\|_{L^2}\Big{)}.\label{gard_pre_final1}
	\end{align}
	Now, we estimate $\|I_h^ku-\mu(\underline{\theta}_h)\|_{L^2}$ using the following dual problem: given $\underline{q}_h:= I_h^ku-\mu(\underline{\theta}_h)$, find $\underline{\phi}_h\in\underline{U}_{h,0}^k$ such that
	\begin{equation}\label{dual_l2_est}
		\mathcal{N}_h^{\text{lin}}(u; \underline{v}_h,\underline{\phi}_h)= (q_h, v_h)\quad \forall\underline{v}_h\in \underline{U}_{h,0}^k.
	\end{equation}
	Choosing $\underline{v}_h= I_h^ku-\mu(\underline{\theta}_h)$ in the above equation, using the definition \eqref{nonlinear_map_mu} of $\mu$, the idea of splitting \eqref{garding_est_1} and \eqref{com_res_est}, we obtain
	\begin{align*}
		&\|I_h^ku-\mu(\underline{\theta}_h)\|_{L^2(\Omega)}^2=\mathcal{N}_h^{\text{lin}}(u; I_h^ku-\mu(\underline{\theta}_h),\underline{\phi}_h)\notag\\ &=\mathcal{N}_h^{\text{lin}}(u;I_h^ku-\underline{\theta}_h, \underline{\phi}_h)+ \mathcal{N}_h(\underline{\theta}_h;\underline{\theta}_h, \underline{\phi}_h)-l(\underline{\phi}_h)\notag\\
		&\leq C\big{(}C_a \big{(}\max_{T\in\cT_h}h_T^{-d/3}\big{)}\left(\|\underline{\theta}_h-I_h^k u\|_{1,h}^2+\|\underline{\theta}_h-I_h^k u\|_{1,h}^3\right)\|\underline{\phi}_h\|_{1,h}\notag\\
		&\quad+h^{r+1}\|\underline{\phi}_h\|_{1,h}+h^{r+1-d/6}\|\underline{\theta}_h-I_h^k u\|_{1,h}\|\underline{\phi}_h\|_{1,h}\big{)}.
	\end{align*}
	Using the a priori bound $\|\underline{\phi}_h\|_{1,h}\lesssim \|I_h^ku-\mu(\underline{\theta}_h)\|_{L^2}$ of \eqref{dual_l2_est}, we obtain
	\begin{align}
		\|I_h^ku-\mu(\underline{\theta}_h)\|_{L^2}&\leq C\Big{(}C_a \big{(}\max_{T\in\cT_h}h_T^{-d/3}\big{)}\left(\|\underline{\theta}_h-I_h^k u\|_{1,h}^2+\|\underline{\theta}_h-I_h^k u\|_{1,h}^3\right)\notag\\
		&\quad\qquad+h^{r+1}+h^{r+1-d/6}\|\underline{\theta}_h-I_h^k u\|_{1,h}\Big{)}.\label{gard_pre_final2}
	\end{align}
	Finally, the above estimations \eqref{gard_pre_final1} and \eqref{gard_pre_final2} lead to
	\begin{align}
		\|I_h^ku-\mu(\underline{\theta}_h)\|_{1,h}
		&\leq \tilde{C}\Big{(}h^{r+1}+h^{r+1-d/6}\|\underline{\theta}_h-I_h^k u\|_{1,h}\notag\\
		&\quad\qquad+(\max_{T\in\cT_h}h_T^{-d/3})\left(\|\underline{\theta}_h-I_h^k u\|_{1,h}^2+\|\underline{\theta}_h-I_h^k u\|_{1,h}^3\right)\Big{)}\label{gard_final_est}
	\end{align}
	for some positive constant $\tilde{C}$ independent of $h$ but  depending on $u$ and $a(u)$.
	Choose $h_{*}$ such that
	\begin{align*}
		(1+2\tilde{C}h_{*}^{r+1-d/6}+4\tilde{C}^2h_{*}^{r+1-d/3}+8\tilde{C}^3h_{*}^{2r+2-d/3})\leq 2.
	\end{align*}
	This implies $(1+2\tilde{C}h^{r+1-d/6}+4\tilde{C}^2h^{r+1-d/3}+8\tilde{C}^3h^{2r+2-d/3})\leq 2$ whenever $h\leq h_{*}$.
	Thus if $\|I_h^ku-\underline{\theta}_h\|_{1,h}\leq R(h):=2\tilde{C}h^{r+1}$, then using Assumption~Q.4 of quasi-uniformity, \eqref{gard_final_est} yields
	\begin{align*}
		&\|I_h^ku-\mu(\underline{\theta}_h)\|_{1,h}\leq \tilde{C}\left(h^{r+1}+2\tilde{C}h^{2r+2-d/6}+4\tilde{C}^2h^{2r+2-d/3}+8\tilde{C}^3h^{3r+3-d/3}\right)\\
		&\quad\leq \tilde{C}h^{r+1}\left(1+2\tilde{C}h^{r+1-d/6}+4\tilde{C}^2h^{r+1-d/3}+8\tilde{C}^3h^{2r+2-d/3}\right)\leq \tilde{C}h^{r+1}\times 2= R(h).
	\end{align*}
	Thus for sufficiently small $h$ $(h\leq h_{*})$, there exist a ball $D(I_h^ku;R(h))$ of radius $R(h)=2\tilde{C}h^{r+1}$ with center at $I_h^ku$ and the following result holds
	$$\|I_h^ku-\underline{\theta}_h\|_{1,h}\leq R(h)\Rightarrow \|I_h^ku-\mu(\underline{\theta}_h)\|_{1,h}\leq R(h).$$
	Hence $\mu$ is a  map from a closed and bounded (compact) convex ball to itself. Therefore, using the Brouwer fixed point theorem, it has a fixed point. This completes the proof.
\end{proof}

To prove the unique fixed point of $\mu$, we show the following contraction result:
\begin{theorem}[Contraction result]\label{thm_cont_mu}
	Let $u\in H^1_0(\Omega)$ be a solution for \eqref{quasilin_pde}. Assume $u\in H^{r+2}(\cT_h)$ for some $r\in\{d-2,\ldots,k\}$. Let $\underline{\theta}_1,\, \underline{\theta}_2\in D(I_h^ku;R(h))$. For sufficiently small $h$, the following contraction result holds:
	$$\|\mu(\underline{\theta}_1)-\mu(\underline{\theta}_2)\|_{1,h}\leq Ch^{r+1-d/3} \|\underline{\theta}_1-\underline{\theta}_2\|_{1,h}.$$
\end{theorem}

\begin{proof}
	Let $\underline{\theta}_1,\, \underline{\theta}_2\in D(I_h^ku;R(h))$, then $\mu(\underline{\theta}_1)$ and $\mu(\underline{\theta}_1)$ satisfy \eqref{nonlinear_map_mu}. That is,
	\begin{align}
		\mathcal{N}_h^{\text{lin}}(u;I_h^ku-\mu(\underline{\theta}_1), \underline{v}_h)= \mathcal{N}_h^{\text{lin}}(u;I_h^ku-\underline{\theta}_1, \underline{v}_h)+ \mathcal{N}_h(\underline{\theta}_1;\underline{\theta}_1, \underline{v}_h)-l(\underline{v}_h),\label{defn_mu_1}\\
		\mathcal{N}_h^{\text{lin}}(u;I_h^ku-\mu(\underline{\theta}_2), \underline{v}_h)= \mathcal{N}_h^{\text{lin}}(u;I_h^ku-\underline{\theta}_2, \underline{v}_h)+ \mathcal{N}_h(\underline{\theta}_2;\underline{\theta}_2, \underline{v}_h)-l(\underline{v}_h).\label{defn_mu_2}
	\end{align}
	Using \GardingsType inequality \eqref{Garding_infsup}, replacing $\underline{\theta}_h$ by $\mu(\underline{\theta}_1)-\mu(\underline{\theta}_2)$ with $\|\underline{v}_h\|_{1,h}=1$, we have
	\begin{align}
		C_1\|\mu(\underline{\theta}_1)-\mu(\underline{\theta}_2)\|_{1,h}\leq \mathcal{N}_h^{\text{lin}}(u;\mu(\underline{\theta}_1)-\mu(\underline{\theta}_2), \underline{v}_h)+ C_2\|\mu(\underline{\theta}_1)-\mu(\underline{\theta}_2)\|_{L^2},\label{mu_garding}
	\end{align}
	where we understand $\|\underline{v}_h\|_{L^2}$ by $\|v_h\|_{L^2(\Omega)}$ in the above term of $\|\mu(\underline{\theta}_1)-\mu(\underline{\theta}_2)\|_{L^2}$.
	From the definition of $\mu$ and subtracting \eqref{defn_mu_1} with \eqref{defn_mu_2}, we get
	\begin{align}
		\mathcal{N}_h^{\text{lin}}(u;\mu(\underline{\theta}_2)-\mu(\underline{\theta}_1), \underline{v}_h)= \mathcal{N}_h^{\text{lin}}(u;\underline{\theta}_2-\underline{\theta}_1, \underline{v}_h)+ \mathcal{N}_h(\underline{\theta}_1;\underline{\theta}_1, \underline{v}_h)- \mathcal{N}_h(\underline{\theta}_2;\underline{\theta}_2, \underline{v}_h).\label{mu_diff}
	\end{align}
	Using the definition $\mathcal{F}_h$ of \eqref{defn_Fh}, the last two terms of the above equation \eqref{mu_diff} yield
	\begin{align}
		&\mathcal{N}_h(\underline{\theta}_1;\underline{\theta}_1, \underline{v}_h)- \mathcal{N}_h(\underline{\theta}_2;\underline{\theta}_2, \underline{v}_h)\notag\\
		&=\langle \mathcal{F}_h(\underline{\theta}_1),\underline{v}_h\rangle-\langle \mathcal{F}_h(\underline{\theta}_2),\underline{v}_h\rangle+s_h(\underline{\theta}_{1}- \underline{\theta}_{2},\underline{v}_h)\notag\\
		&=\left(\langle \mathcal{F}_h(\underline{\theta}_1),\underline{v}_h\rangle-\langle \mathcal{F}_h(I_h^k u),\underline{v}_h\rangle+s_h(\underline{\theta}_{1}- I_h^k u,\underline{v}_h)\right)\notag\\
		&\quad-\left(\langle \mathcal{F}_h(\underline{\theta}_2),\underline{v}_h\rangle-\langle \mathcal{F}_h(I_h^k u),\underline{v}_h\rangle+s_h( \underline{\theta}_{2}-I_h^k u,\underline{v}_h)\right).
		\label{Gh_diff_cont}
	\end{align}
	Set $\underline{\chi}_{1h}:=\underline{\theta}_1-I_h^k u,  \underline{\chi}_{2h}:=\underline{\theta}_2-I_h^k u$ and $\check{u}_T:=R_T^{k+1}I_T^ku$. From the expansion of $\mathcal{F}_h$ of \eqref{Gh_expan}, we rewrite the above equation \eqref{Gh_diff_cont} as
	\begin{align}
		&\mathcal{N}_h(\underline{\theta}_1;\underline{\theta}_1, \underline{v}_h)- \mathcal{N}_h(\underline{\theta}_2;\underline{\theta}_2, \underline{v}_h)
		-\widetilde{\mathcal{N}}_h^{\text{lin}}(I_h^k u;\underline{\theta}_1-\underline{\theta}_2,\underline{v}_h)\notag\\
		&=\Bigg{(}\sit\tilde{a}_{uu}(\check{u}_T)(R_T^{k+1}\underline{\chi}_{1T})^2\nabla \check{u}_T{\cdot} \Gtk\underline{v}_T\dx\notag\\
		&\qquad\qquad\qquad-\sit\tilde{a}_{uu}(\check{u}_T)(R_T^{k+1}\underline{\chi}_{2T})^2\nabla \check{u}_T{\cdot} \Gtk\underline{v}_T\dx\Bigg{)}\notag\\
		&\quad+\Bigg{(}\sit a_u(\check{u}_T)R_T^{k+1}\underline{\chi}_{1T}\Gtk\underline{\chi}_{1T}{\cdot}\Gtk\underline{v}_T\dx\notag\\
		&\qquad\qquad\qquad-\sit a_u(\check{u}_T)R_T^{k+1}\underline{\chi}_{2T}\Gtk\underline{\chi}_{2T}{\cdot}\Gtk\underline{v}_T\dx\Bigg{)}\notag\\
		&\qquad+ \Bigg{(}\sit \tilde{a}_{uu}(\check{u}_T)(R_T^{k+1}\underline{\chi}_{1T})^2\Gtk\underline{\chi}_{1T}{\cdot} \Gtk\underline{v}_T\dx\notag\\
		&\qquad\qquad-\sit \tilde{a}_{uu}(\check{u}_T)(R_T^{k+1}\underline{\chi}_{2T})^2\Gtk\underline{\chi}_{2T}{\cdot} \Gtk\underline{v}_T\dx\Bigg{)}=:T_1+T_2+T_3.\label{Bh_diff_lin}
	\end{align}
	To estimate the above terms, we use the identities
	\begin{align}
		&a^2-b^2=(a-b)(a+b),\quad a_1b_1-a_2b_2= a_1(b_1-b_2)+ (a_1-a_2)b_2,\label{sim_iden_1}\\
		&a_1^2b_1-a_2^2b_2= a_1^2(b_1-b_2)+ (a_1-a_2)(a_1+a_2)b_2.\label{sim_iden_2}
	\end{align}
	Using the above identity \eqref{sim_iden_1} and estimates similar to Lemma~\ref{lem_imp_ineq}(a), the terms $T_1$ and $T_2$ of \eqref{Bh_diff_lin} yield
	\begin{align*}
		&\sit \tilde{a}_{uu}(\check{u}_T) (R_T^{k+1}\underline{\theta}_{1T}-R_T^{k+1}\underline{\theta}_{2T})(R_T^{k+1}\underline{\theta}_{1T}-\check{u}_T+ R_T^{k+1}\underline{\theta}_{2T}-\check{u}_T)\nabla \check{u}_T{\cdot} \Gtk\underline{v}_T\dx\\
		&\quad+ \sit a_{u}(\check{u}_T)(\check{u}_T-R_T^{k+1}\underline{\theta}_{1T})\Gtk(\underline{\theta}_{1T}-\underline{\theta}_{2T}){\cdot} \Gtk\underline{v}_T\dx\notag\\
		&\qquad+ \sit a_{u}(\check{u}_T)(R_T^{k+1}\underline{\theta}_{1T}-R_T^{k+1}\underline{\theta}_{2T})\Gtk(I_T^ku-\underline{\theta}_{2T}){\cdot} \Gtk\underline{v}_T\dx\\
		&\leq C_a C\big{(}\max_{T\in\cT_h}h_T^{-d/6}\big{)}\|\underline{\theta}_{1}- \underline{\theta}_{2}\|_{1,h} \left(\|I_h^ku-\underline{\theta}_{2}\|_{1,h}+\|I_h^ku-\underline{\theta}_{1}\|_{1,h}\right) \|\underline{v}_h\|_{1,h}.
	\end{align*}
	The above identity \eqref{sim_iden_2} and estimates similar to Lemma~\ref{lem_imp_ineq} lead to an estimate for $T_3$ term of \eqref{Bh_diff_lin} as
	\begin{align*}
		T_3&=\sit \tilde{a}_{uu}(\check{u}_T)\Big{(} (\check{u}_T-R_T^{k+1}\underline{\theta}_{2T})^2\Gtk (\underline{\theta}_{1T}-\underline{\theta}_{2T})+(R_T^{k+1}\underline{\theta}_{1T}-R_T^{k+1}\underline{\theta}_{2T})\times\\
		&\qquad(\check{u}_T-R_T^{k+1}\underline{\theta}_{2T}+ \check{u}_T-R_T^{k+1}\underline{\theta}_{1T}) \Gtk(I_T^ku-\underline{\theta}_{1T})\Big{)}{\cdot} \Gtk\underline{v}_T\dx\\
		&\leq C_a C\big{(}\max_{T\in\cT_h}h_T^{-d/3}\big{)}\|\underline{\theta}_{1}- \underline{\theta}_{2}\|_{1,h} \left(\|I_h^ku-\underline{\theta}_{2}\|_{1,h}^2+\|I_h^ku-\underline{\theta}_{1}\|_{1,h}^2\right) \|\underline{v}_h\|_{1,h}.
	\end{align*}
	Combine the estimation for \eqref{Bh_diff_lin} in \eqref{mu_diff} with $\|\underline{v}_h\|_{1,h}=1$, it yields
	\begin{align}
		&\mathcal{N}_h^{\text{lin}}(u;\mu(\underline{\theta}_2)-\mu(\underline{\theta}_1), \underline{v}_h)=\left(\mathcal{N}_h^{\text{lin}}(u;\underline{\theta}_{2}- \underline{\theta}_{1},\underline{v}_h)-\widetilde{\mathcal{N}}_h^{\text{lin}}(I_h^ku;\underline{\theta}_{2}- \underline{\theta}_{1}, \underline{v}_h)\right)\notag\\
		&\quad+\left(\mathcal{N}_h(\underline{\theta}_1;\underline{\theta}_1, \underline{v}_h)- \mathcal{N}_h(\underline{\theta}_2;\underline{\theta}_2, \underline{v}_h)- \widetilde{\mathcal{N}}_h^{\text{lin}}(I_h^ku;\underline{\theta}_{1}- \underline{\theta}_{2}, \underline{v}_h)\right)\notag\\
		&\leq Ch^{r+1-d/6}\|\underline{\theta}_{1}- \underline{\theta}_{2}\|_{1,h}+C_a C\big{(}\max_{T\in\cT_h}h_T^{-d/3}\big{)}\|\underline{\theta}_{1}- \underline{\theta}_{2}\|_{1,h} \bigg{(}\|I_h^ku-\underline{\theta}_{1}\|_{1,h}\notag\\
		&\hspace*{3cm}+\|I_h^ku-\underline{\theta}_{1}\|_{1,h}^2
		+\|I_h^ku-\underline{\theta}_{2}\|_{1,h}+\|I_h^ku-\underline{\theta}_{2}\|_{1,h}^2\bigg{)}.\label{cts_lin_nonselfadj05}
	\end{align}
	To obtain the estimation for $L^2$-term $\|\mu(\underline{\theta}_1)-\mu(\underline{\theta}_2)\|_{L^2}$ of \eqref{mu_garding}, we consider the dual linear problem: given $\underline{q}_h:=\mu(\underline{\theta}_1)-\mu(\underline{\theta}_2)$, find $\underline{\phi}_h\in \underline{U}_{h,0}^k$ such that
	\begin{align}
		\mathcal{N}_h^{\text{lin}}(u; \underline{v}_h,\underline{\phi}_h)= (q_h,v_h)\quad \forall \underline{v}_h\in \underline{U}_{h,0}^k.\label{dual_cont_prob} 
	\end{align}
	Choose $\underline{v}_h= \mu(\underline{\theta}_1)-\mu(\underline{\theta}_2)$ in the above equation \eqref{dual_cont_prob} and use \eqref{cts_lin_nonselfadj05} to obtain
	\begin{align*}
		&\|\mu(\underline{\theta}_1)-\mu(\underline{\theta}_2)\|_{L^2}^2= \mathcal{N}_h^{\text{lin}}(u; \mu(\underline{\theta}_1)-\mu(\underline{\theta}_2),\underline{\phi}_h)\\
		&\leq Ch^{r+1-d/6}\|\underline{\theta}_{1}- \underline{\theta}_{2}\|_{1,h}\|\underline{\phi}_h\|_{1,h}+C_a C\big{(}\max_{T\in\cT_h}h_T^{-d/3}\big{)}\|\underline{\theta}_{1}- \underline{\theta}_{2}\|_{1,h}\\ &\quad\times\Big{(}\|I_h^ku-\underline{\theta}_{1}\|_{1,h}+\|I_h^ku-\underline{\theta}_{1}\|_{1,h}^2+\|I_h^ku-\underline{\theta}_{2}\|_{1,h}+\|I_h^ku-\underline{\theta}_{2}\|_{1,h}^2\Big{)}\|\underline{\phi}_h\|_{1,h}.
	\end{align*}
	The a priori bound $\|\underline{\phi}_h\|_{1,h}\leq C \|\mu(\underline{\theta}_1)-\mu(\underline{\theta}_2)\|_{L^2}$ of \eqref{dual_cont_prob} leads to
	\begin{align}
		&\|\mu(\underline{\theta}_1)-\mu(\underline{\theta}_2)\|_{L^2}\leq Ch^{r+1-d/6}\|\underline{\theta}_{1}- \underline{\theta}_{2}\|_{1,h}+ C_a C\|\underline{\theta}_{1}- \underline{\theta}_{2}\|_{1,h}\big{(}\max_{T\in\cT_h}h_T^{-d/3}\big{)}\notag\\ &\quad\times\Big{(}\|I_h^ku-\underline{\theta}_{1}\|_{1,h}+\|I_h^ku-\underline{\theta}_{1}\|_{1,h}^2+\|I_h^ku-\underline{\theta}_{2}\|_{1,h}+\|I_h^ku-\underline{\theta}_{2}\|_{1,h}^2\Big{)}.\label{cts_lin_nonselfadj08}
	\end{align}
	Using \eqref{cts_lin_nonselfadj05} and \eqref{cts_lin_nonselfadj08}, we obtain from \eqref{mu_garding} as
	\begin{align*}
		&\|\mu(\underline{\theta}_1)-\mu(\underline{\theta}_2)\|_{1,h}\leq Ch^{r+1-d/6}\|\underline{\theta}_{1}- \underline{\theta}_{2}\|_{1,h}+C_a C\|\underline{\theta}_{1}- \underline{\theta}_{2}\|_{1,h}\big{(}\max_{T\in\cT_h}h_T^{-d/3}\big{)}\notag\\ &\quad\times\Big{(}\|I_h^ku-\underline{\theta}_{1}\|_{1,h}+\|I_h^ku-\underline{\theta}_{1}\|_{1,h}^2+\|I_h^ku-\underline{\theta}_{2}\|_{1,h}+\|I_h^ku-\underline{\theta}_{2}\|_{1,h}^2\Big{)}.
	\end{align*}
	Since $\underline{\theta}_{1},\, \underline{\theta}_{2}\in D(I_h^ku;R(h))$ with $R(h)=2\tilde{C}h^{r+1}$, we have
	\begin{equation*}
		\|I_h^ku-\underline{\theta}_{1}\|_{1,h}\leq 2\tilde{C}h^{r+1}\quad \text{and}\quad \|I_h^ku-\underline{\theta}_{2}\|_{1,h}\leq 2\tilde{C}h^{r+1}.
	\end{equation*}
	This implies
	\begin{equation*}
		\|\mu(\underline{\theta}_1)-\mu(\underline{\theta}_2)\|_{1,h}\leq Ch^{r+1-d/3}\|\underline{\theta}_{1}- \underline{\theta}_{2}\|_{1,h}    
	\end{equation*}
	for sufficiently small mesh parameter $h$.
	This completes the proof.
\end{proof}

For sufficiently small $h$, Theorem~\ref{thm_cont_mu} proves the local uniqueness of the fixed point of $\mu$, and hence the local uniqueness of the solution of \eqref{hho_dis_quasi}.

\noindent{\bf Error estimate for quasilinear problem}.
Adding and subtracting $\Ghk I_h^k u$, using triangle inequality and Theorem~\ref{thm_fixed_point}, we have the following error estimate:
\begin{theorem}[Error estimate]\label{thm_err_est_quasi}
	Let $u\in H^1_0(\Omega)$ be the solution of nonlinear problem \eqref{quasilin_pde} and $\underline{u}_h\in\underline{U}_{h,0}^k$ be the solution of the discrete problem \eqref{hho_dis_quasi}. Assume $u\in H^{r+2}(\cT_h)$ for some $r\in\{d-2,\ldots,k\}$. Then for sufficiently small $h$, we have
	\begin{align}\label{err_est_quasi}
		\|\nabla u-\Ghk\underline{u}_h\|\leq Ch^{r+1}.
	\end{align}
\end{theorem}


\section{Numerical experiments}\label{Sec:Num_Quasilin}
\begin{figure}
	\begin{center}
		\subfloat[]{\includegraphics[width=0.25\textwidth]{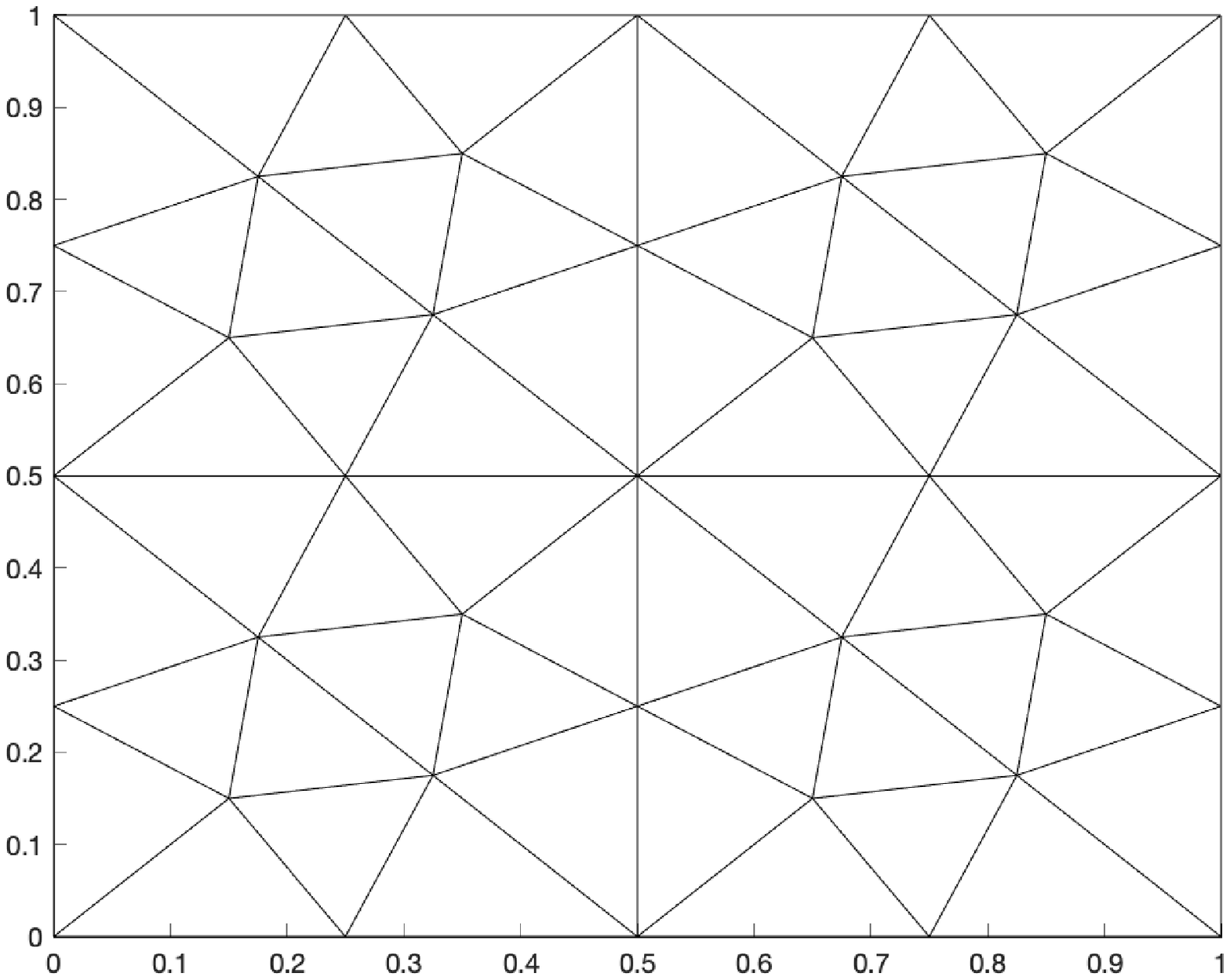}}
		\subfloat[]{\includegraphics[width=0.25\textwidth]{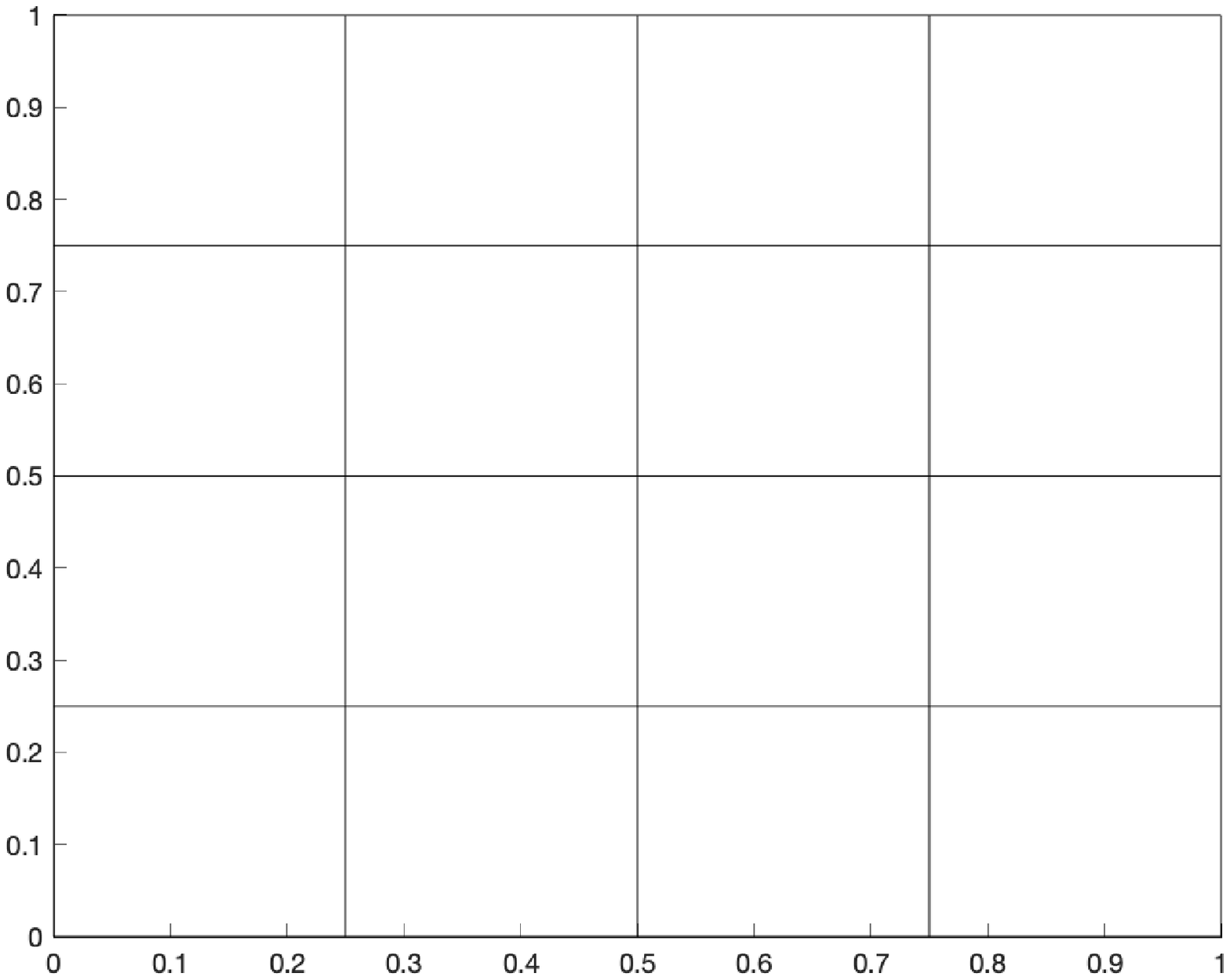}}
		\subfloat[]{\includegraphics[width=0.25\textwidth]{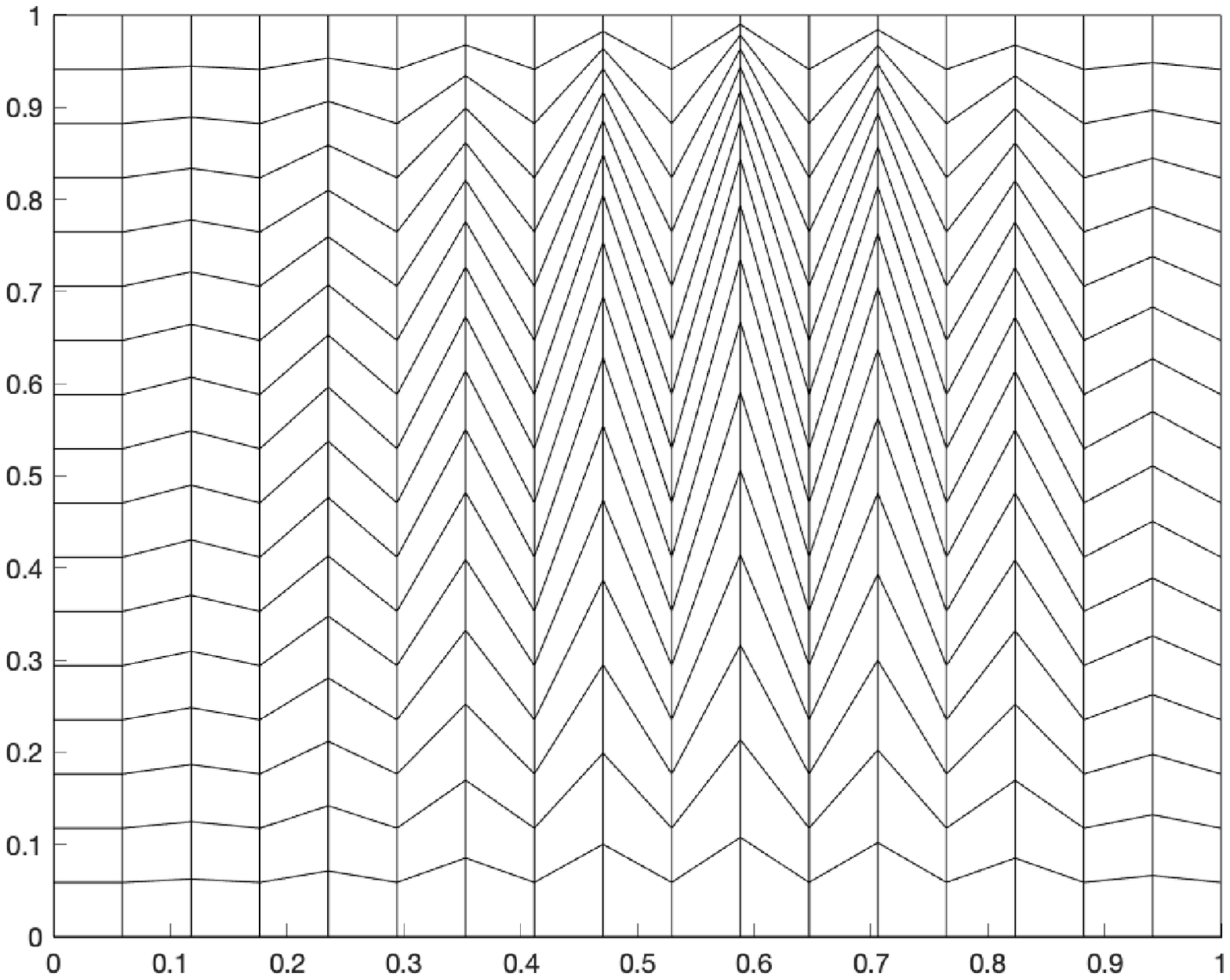}}
		\subfloat[]{\includegraphics[width=0.25\textwidth]{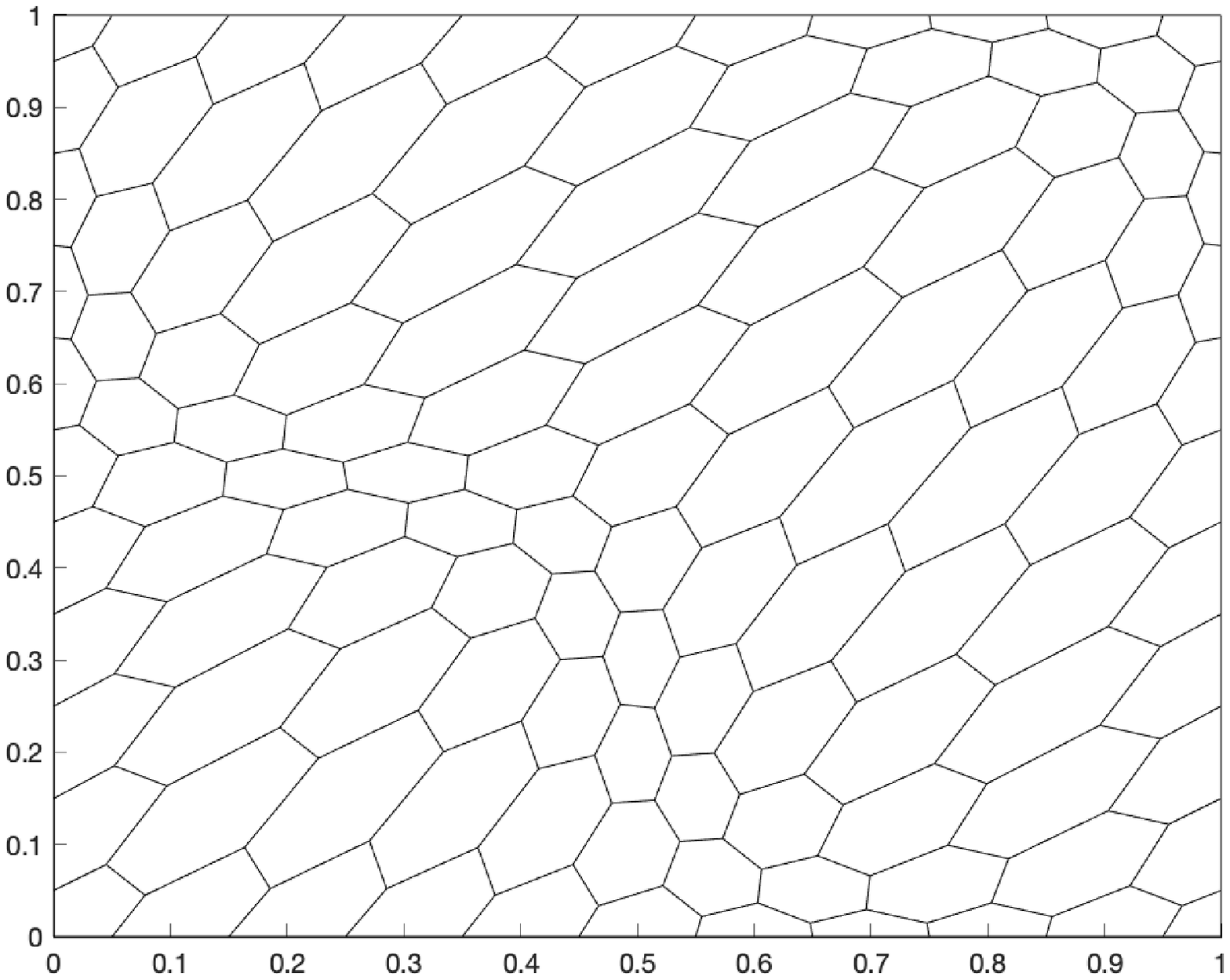}}
		\caption{(a) Triangular, (b)  Cartesian, (c) Kershaw, and (d) hexagonal initial meshes.}
		\label{fig:Triag_Cart_Meshes}
	\end{center}
\end{figure}

\begin{figure}
	\begin{center}
		\subfloat[]{\includegraphics[height=0.4\textwidth,width=0.5\textwidth]{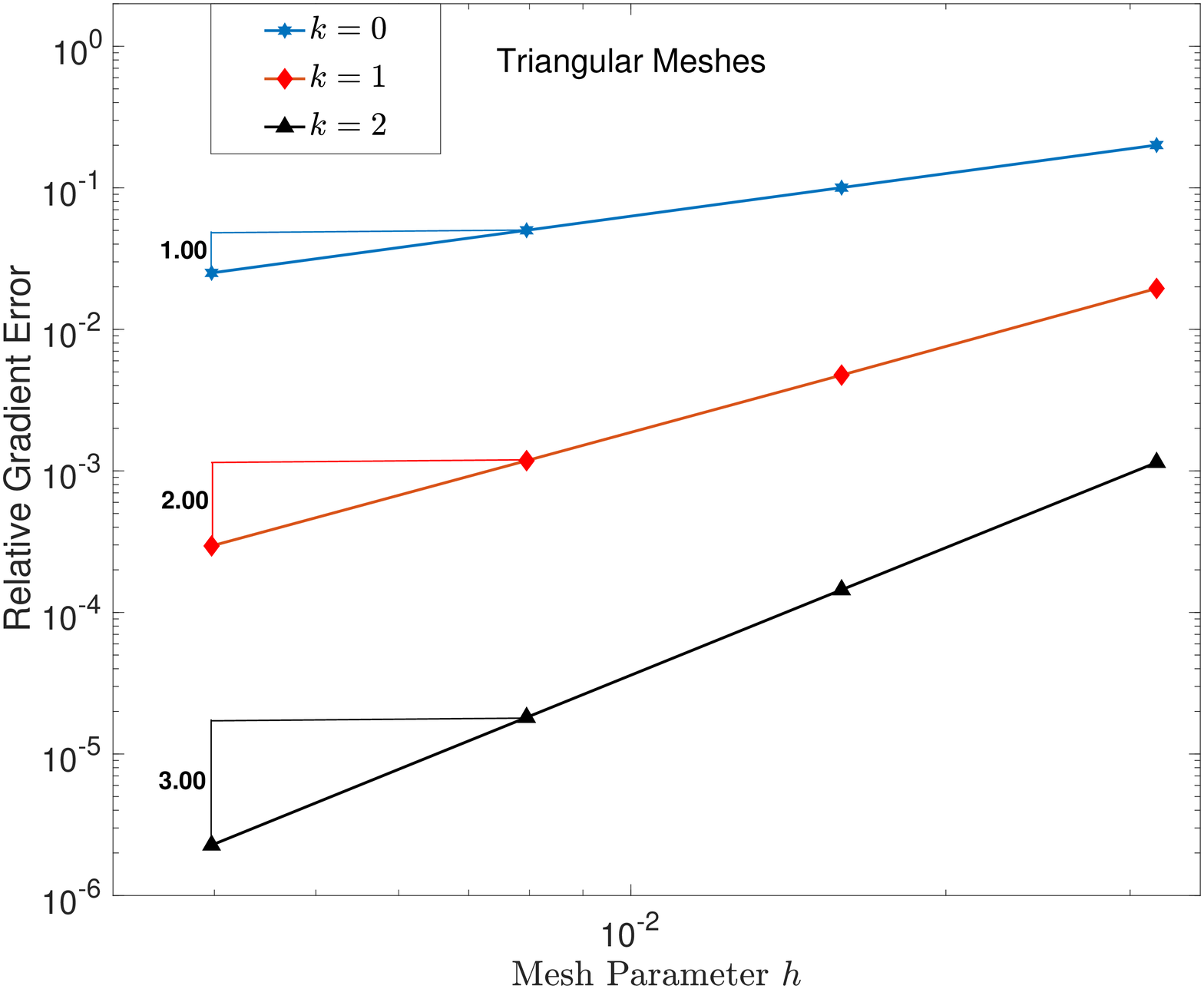}}
		\subfloat[]{\includegraphics[height=0.4\textwidth,width=0.5\textwidth]{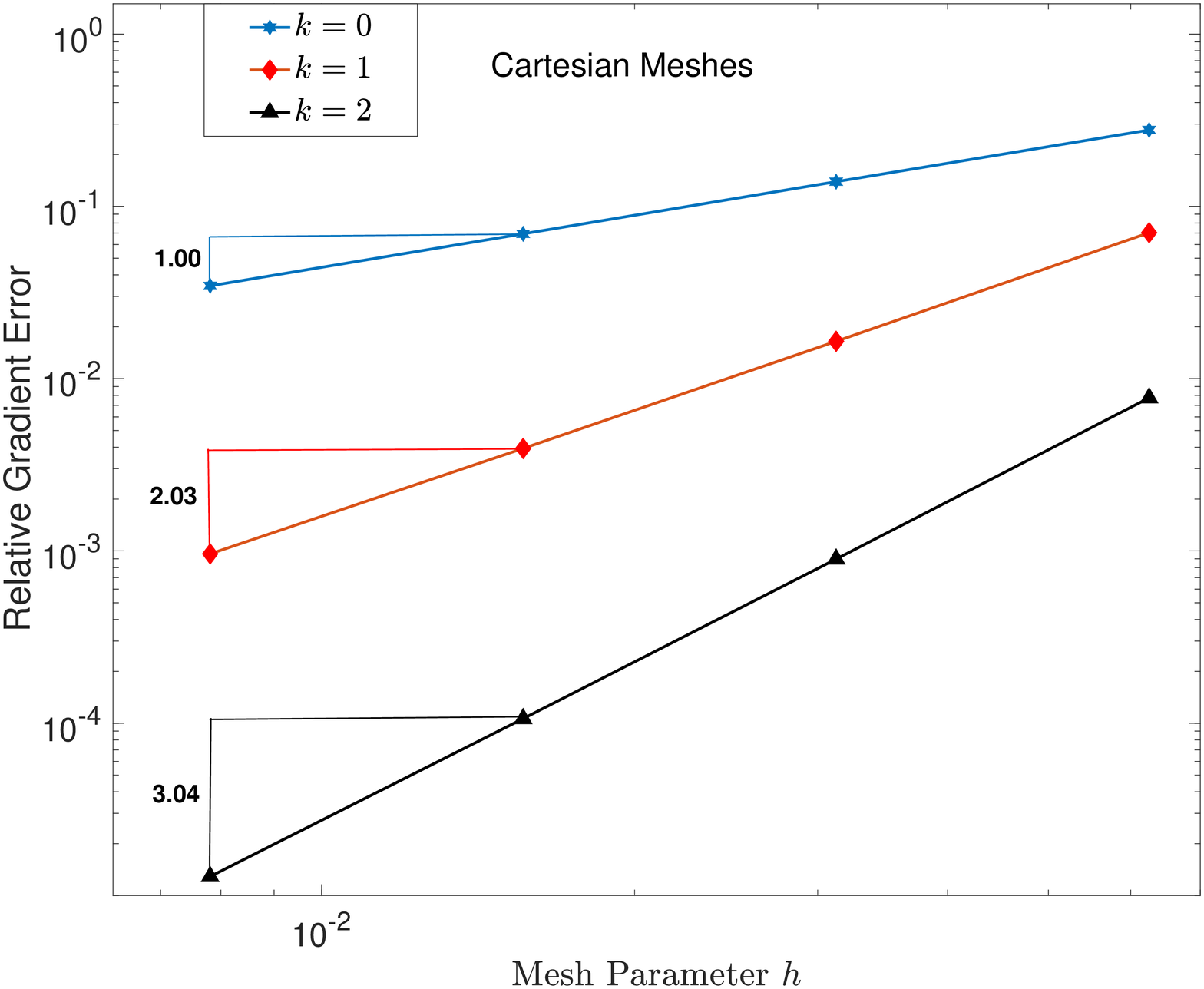}}\\
		\subfloat[]{\includegraphics[height=0.4\textwidth,width=0.5\textwidth]{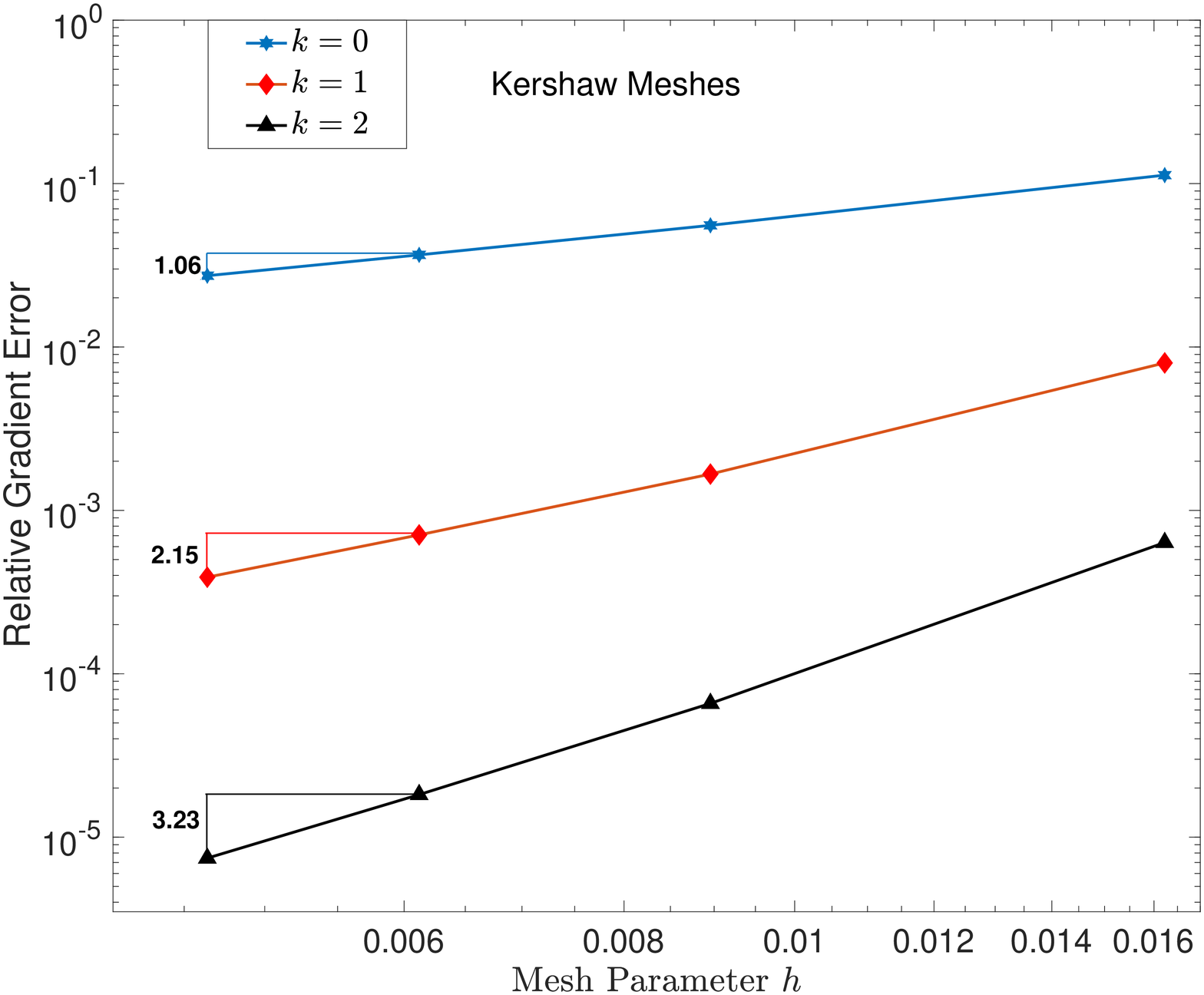}}
		\subfloat[]{\includegraphics[height=0.4\textwidth,width=0.5\textwidth]{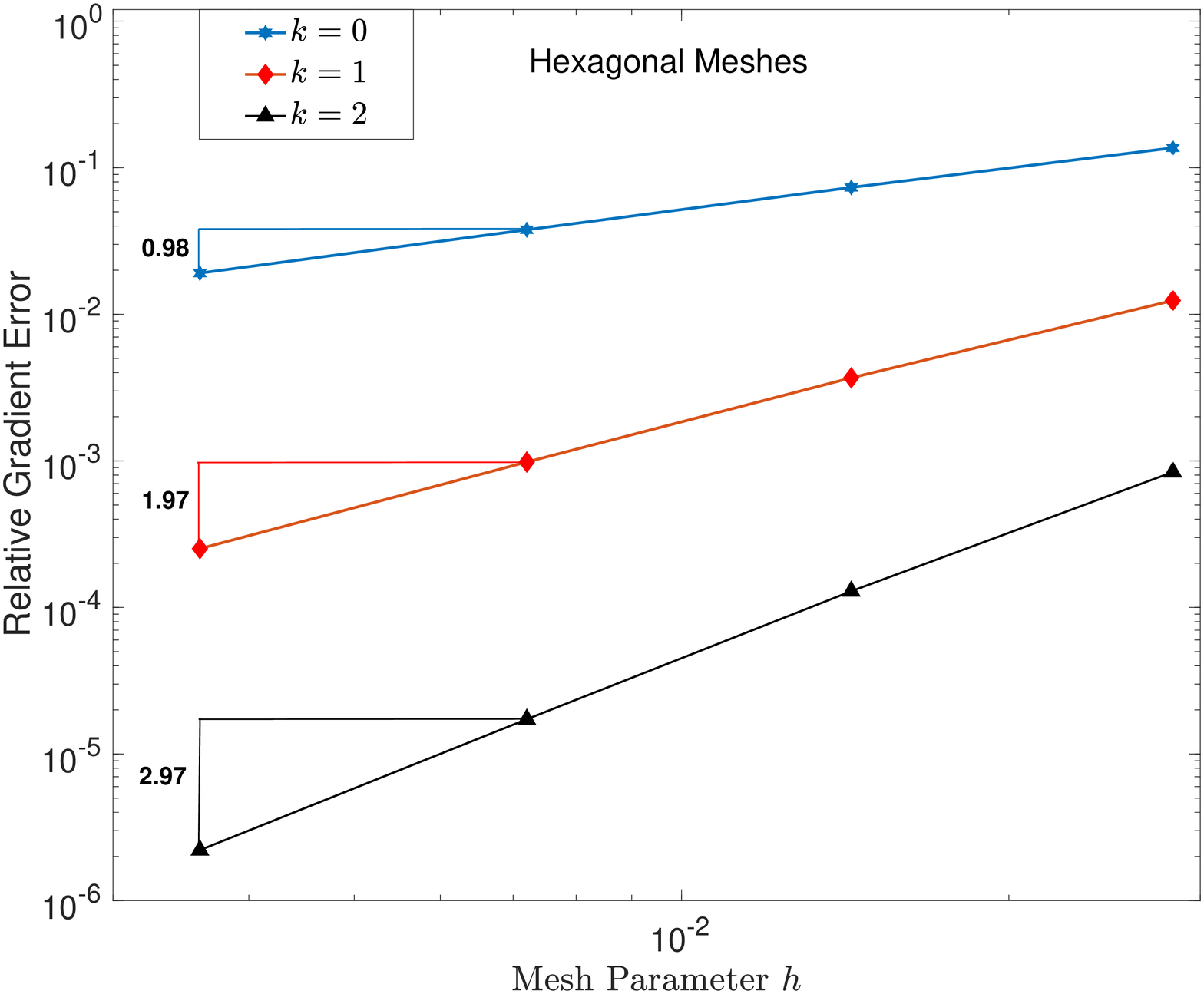}}
		\caption{Convergence histories for the relative gradient error on (a) Triangular, (b) Cartesian   (c) Kershaw and (d) hexagonal meshes.}
		\label{fig:Quasi_Conv_His_Triag_Cart}
	\end{center}
\end{figure}

In this section, numerical experiments are performed using the hybrid high-order approximation \eqref{hho_dis_quasi} for the quasilinear problem \eqref{quasilin_pde}. We consider the following nonlinear problem \cite{Gudi_AKP_07_DG_quasi}:
\begin{subequations}
	\begin{align*}
		-\nabla{\cdot}\left((1+u)\nabla u\right) &= f\quad\text{in}~~ \Omega,\\
		u&=0\quad\text{on}~~ \partial\Omega,
	\end{align*}
\end{subequations}
where $\Omega:=(0,1)\times (0,1)\subset \bR^2$, and the source term $f$ is taken in such a way that the exact solution reads $u(x,y)=x(1-x)y(1-y)$. We rewrite the nonlinear map \eqref{nonlinear_map_mu} in order to obtain practical iterative solutions. In the computation, we do not demand the exact solution $u$, and this is replaced by the previous step computed (initial guess) solution. We start with an initial guess $u_h^0\in \underline{U}_{h,0}^k$ obtained from solving the Dirichlet Poisson problem $-\Delta u=f$ with the same load function $f$ as defined above. The $(n+1)$-th iteration is given by, for all $\underline{v}_h \in\underline{U}_{h,0}^k$,
\begin{align*}
	\widetilde{\mathcal{N}}_h^{\text{lin}}(\underline{u}_h^{n};\underline{u}_h^{n+1},\underline{v}_h)=\widetilde{\mathcal{N}}_h^{\text{lin}}(\underline{u}_h^{n};\underline{u}_h^{n},\underline{v}_h)-\mathcal{N}_h(\underline{u}_h^n;\underline{u}_h^n,\underline{v}_h)+l(\underline{v}_h)\; \text{ for }n=0,1,2,\ldots,
\end{align*} 
where the linearized  $\widetilde{\mathcal{N}}_h^{\text{lin}}$ and nonlinear $\cN_h$  forms are as defined in \eqref{defn_Bh_lin} and \eqref{defn_Bh} respectively. The stopping criteria is prescribed by a tolerance $10^{-10}$ for the relative gradient error for the difference of two successive iterative solutions as $\|\Ghk(\underline{u}_h^{n+1}-\underline{u}_h^{n})\|/\|\Ghk\underline{u}_h^{n+1}\|\leq 10^{-10}$. We consider the triangular, Cartesian, Kershaw and hexagonal mesh families for numerical experiments  which are depicted in Figure~\ref{fig:Triag_Cart_Meshes}. The triangular, Cartesian and Kershaw mesh families are discussed in \cite{Her_Hub_FVCA_mesh_08}, and hexagonal mesh family is introduced in \cite{Piet_Lem_HexaMesh_15}. 
The experiments are performed in Matlab. Some of the basic implementation methodologies for hybrid high-order method are adopted from \cite{Piet_Jero_HHO_Book_20,Cicu_Piet_Ern_18_Comp_HHO,Piet_Ern_Lem_14_arb_local}.
It has been observed that the iterative step terminates within $4$-steps using the above stopping criteria.
The experimental rate of convergence is computed as
\begin{align*}
	\texttt{rate}(\ell):=\log \big(e_{h_{\ell}}/e_{h_{\ell-1}} \big)/\log \big(h_\ell/h_{\ell-1} \big)\text{ for } \ell=1,2,3,\ldots,
\end{align*}
where $e_{h_{\ell}}$ and $e_{h_{\ell-1}}$ are the errors associated to two consecutive meshsizes $h_\ell$ and $h_{\ell-1}$ respectively.
In Figure~\ref{fig:Quasi_Conv_His_Triag_Cart}, we have plotted the convergence histories for  the relative reconstructed gradient error $e_h=\|\nabla u-\Ghk\underline{u}_h\|/\|\nabla u\|$ as a function of mesh parameter $h$ on the sequence of triangular, Cartesian, Kershaw and  hexagonal meshes for the polynomial degree $k=0,1,2$. 
The convergence rate for the polynomial degree $k=0,1,2$ are respectively close to $1,2,3$ for each of the mesh. The convergence rates are in line with the theoretical convergence found in Theorem~\ref{thm_err_est_quasi}.

\section{Conclusions}\label{sec:Conclusion}
In this paper, we have discussed a hybrid high-order approximation for the second-order quasilinear elliptic problem of nonmonotone type defined on polytopal domain in $\bR^d, d=2,3$. First, we have deduced existence, uniqueness and error estimate for the discrete solution of a general second-order nonselfadjoint problem. This has helped us to construct a nonlinear map which satisfies contraction property over a small ball. The discrete solution of nonlinear problem is essentially a fixed point of the nonlinear map. The existence, local uniqueness and error estimate for HHO approximation of nonlinear problem are established. The analysis does not require any user specified large penalty parameter unlike the discontinuous Galerkin method of \cite{Gudi_AKP_07_DG_quasi}.  
The analysis also supports lowest-order ($k=0$ when $d=2$) polynomial approximation with linear order of convergence.
It is possible to extend our analysis without much difficulties to the more general nonlinear problem of the type $\nabla{\cdot}(a(x,u)\nabla u)+f(x,u)=0$ with $f\in C^2_b(\bar{\Omega}\times \mathbb{R})$, and when $a(x,u)$ is uniformly bounded positive-definite matrix.


\bigskip
\noindent{\bf Funding}\\
The first author acknowledges the financial support of DST MATRICS grant. The second author acknowledges the financial support of National Board for Higher Mathematics (NBHM) research grant 0204/58/2018/R\&D-II/14746 and DST C.V. Raman grant R(IA)/CVR-PDF/2020.

\bibliography{HHO_Bib}
\bibliographystyle{abbrvnat}

\end{document}